\documentclass[oneside,12pt]{amsart}

\usepackage{amsfonts}
\usepackage{amssymb}
\usepackage{amsxtra}
\usepackage{amstext}
\usepackage[english]{babel}

\usepackage[T1]{fontenc}
\usepackage{latexsym}
\usepackage[latin1]{inputenc}

\newtheorem{theorem}{Theorem}[section]
\newtheorem{lemma}[theorem]{Lemma}
\newtheorem{corollary}[theorem]{Corollary}
\newtheorem{proposition}[theorem]{Proposition}
\newtheorem{example}[theorem]{Example}
\newtheorem{remark}[theorem]{Remark}

\newtheorem{definition}[theorem]{Definition}
\def\bit{\begin{itemize}}
\def\eit{\end{itemize}}
\reversemarginpar   
\def\bc{\begin{center}}
\def\ec{\end{center}}
\def\bthm{\begin{theorem}}
\def\ethm{\end{theorem}}
\def\bcor{\begin{corollary}}
\def\ecor{\end{corollary}}
\def\bprop{\begin{proposition}}
\def\eprop{\end{proposition}}
\def\blem{\begin{lemma}}
\def\elem{\end{lemma}}

\def\brem{\begin{remark}}
\def\erem{\end{remark}}
\def\prf{\noindent{\bf Proof~: }}
\def\bdes{\begin{description}}
\def\edes{\end{description}}
\def\ita{\item[(a)]}
\def\itb{\item[(b)]}
\def\itc{\item[(c)]}

\def\iti{\item[(i)]}
\def\itii{\item[(ii)]}
\def\itiii{\item[(iii)]}

\def\beq{\begin{equation}}
\def\eeq{\end{equation}}

\def\ben{\begin{enumerate}}
\def\een{\end{enumerate}}

\def\beqar{\begin{eqnarray}}
\def\eeqar{\end{eqnarray}}
\def\beqarr{\begin{eqnarray*}}
\def\eeqarr{\end{eqnarray*}}

\def\RR{{\mathbb R}}  

\def\PP{{\mathbb P}}
\def\QQ{{\mathbb Q}}
\def\cA{\mathcal{A}} \def\cB{\mathcal{B}} 
 \def\cE{\mathcal{E}} \def\cF{\mathcal{F}}

\def\cM{\mathcal{M}}  
\def\cP{\mathcal{P}}

\def\wt{\widetilde}

\def\P{{\mathsf P}} 

\def\E{{\mathsf E}} 

\def\eps{\epsilon}

\def\p{\varphi}

\def\part{\partial}

\def\d#1dt{\frac{d#1}{dt}}    

\begin{document}
\title{Three examples of Brownian flows on $\RR$}
\author{Yves Le Jan}
\address{Laboratoire Mod\'elisation stochastique et statistique, Universit\'e Paris-Sud, B\^atiment 425, 91405 Orsay Cedex}
\email{yves.lejan@math.u-psud.fr}
\author{Olivier Raimond}
\address{Laboratoire Modal'X, Universit\'e Paris Ouest Nanterre La D\'efense, B\^atiment G, 200 avenue
de la R\'epublique 92000 Nanterre, France.}
\email{olivier.raimond@u-paris10.fr}

\maketitle

\section*{Abstract}
We show that the only flow solving the stochastic differential equation (SDE) on $\RR$
$$dX_t = 1_{\{X_t>0\}}W_+(dt) + 1_{\{X_t<0\}}dW_-(dt),$$
where $W^+$ and $W^-$ are two independent  white noises, is a coalescing flow we will denote $\p^{\pm}$. The flow $\p^\pm$ is a Wiener solution. Moreover, $K^+=\E[\delta_{\p^\pm}|W_+]$ is the unique solution (it is also a Wiener solution) of the SDE 
$$K^+_{s,t}f(x)=f(x)+\int_s^t K_{s,u}(1_{\RR^+}f')(x)W_+(du)+\frac {1}{2} \int_s^t K_{s,u}f''(x) du$$
for $s<t$, $x\in\RR$ and $f$ a twice continuously differentiable function. A third flow $\p^+$ can be constructed out of the $n$-point motions of $K^+$. This flow is coalescing and its $n$-point motion is given by the $n$-point motions of $K^+$ up to the first coalescing time, with the condition that when two points meet, they stay together. We note finally that $K^+=\E[\delta_{\p^+}|W_+]$.

\section{Introduction}
Our purpose in this paper is to study two very simple one
dimensional SDE's which can be completely solved, although they do
not satisfy the usual criteria. This study is done in the framework
of stochastic flows exposed in \cite{LJR2002,LJR2004a,LJR2005} and \cite{Tsirelson04} (see also \cite{rozov04} and \cite{rozov05} in an SPDE setting). 
There is still a lot to do to understand the nature of these flows, even if one consider only Brownian flows, i.e. flows whose one-point motion is a Brownian motion.
As in our previous study of
Tanaka's equation  (see \cite{LJR2006} and also \cite{hajri11,hajri12}), it is focused on the case where the
singularity is located at an interface between two half lines. It
should be generalizable to various situations. The first SDE
represents the motion of particles driven by two independent white noises
$W_+$ and $W_-$. All particles on the positive half-line are driven
by $W_+$ and therefore move parallel until they hit $0$. $W_-$
drives in the same way the particles while they are on the negative
side. What should happen at the origin is a priori not clear, but we
should assume particles do not spend a positive measure of time there.
The SDE can therefore be written
$$dX_t= 1_{\{X_t>0\}} W_+(dt) + 1_{\{X_t<0\}} W_-(dt)$$ 
and will be shown to have a strong, i.e. $ \sigma(W_+, W_-)$ (Wiener) measurable, solution which is a coalescing flow of continuous maps.
This is the only solution in any reasonable sense, even if one
allows an extension of the probability space, i.e. other sources of
randomness than the driving noises, to define the flow.

If we compare this result with the one obtained in \cite{LJR2006}  for
Tanaka's equation, in which the construction of the coalescing flow requires an additional countable family of independent Bernoulli variable attached to local minima of the noise $W$,  the Wiener measurability may seem somewhat surprising. 
A possible intuitive interpretation is the following: In
the case of Tanaka's equation, if we consider the image of zero, a
choice has to be made at the beginning of every excursion outside of zero of the
driving reflected Brownian motion. In the case of our
SDE, an analogous role is played by excursions of $W^+$ and $W^-$. 
These excursions have to be taken at various levels different of $0$, but the essential point is that at given levels they a.s. never start at the same time. 
And if we could a priori neglect the effect of the excursions of height smaller than some positive
$\epsilon$, the motion of a particle starting at zero would be perfectly determined.

The second SDE is a transport equation, which cannot be induced by a
flow of maps. The matter is dispersed according to the heat equation
on the negative half line and is driven by $W_+$ on the positive
half line. A solution is easily constructed by integrating out $W_-$
in the solution of our first equation. We will prove that also in
this case, there is no other solution, even on an extended
probability space.

The third flow is not related to an SDE. 
It is constructed in a similar way as Arratia flow (see \cite{LJR2004a,arratia,Tsirelson04,Fontes02,Fontes04}) is constructed: 
the $n$-point motion is given by independent Brownian motions that coalesce when they meet (without this condition the matter is dispersed according to the heat equation on the line). 
Using the same procedure, each particle are driven on the negative half line by an independent white noises and on the positive half line by $W_+$. 
This procedure permits to define a coalescing flow of maps, which is not Wiener measurable.

\pagebreak
\section{Notation, definitions and results.}
\subsection{Notation}
\begin{itemize}
\item For $n\ge 1$, $C(\RR^+:\RR^n)$ (resp. $C_b(\RR^+:\RR^n)$) denotes the space of continuous (resp. bounded continuous) functions $f:\RR^+\to \RR^n$.
\item For $n\ge 1$, $C_0(\RR^n)$ is the space of continuous functions $f:\RR^n\to\RR$ converging to $0$ at infinity. It is equipped with the norm $\|f\|_\infty=\sup_{x\in\RR^n}|f(x)|$.
\item For $n\ge 1$, $C^2_0(\RR^n)$ is the space of twice continuously differentiable functions $f:\RR^n\to\RR$  converging to $0$ at infinity as well as their derivatives. It is equipped with the norm $\|f\|_{2,\infty}=\|f\|_\infty + \sum_i\|\partial_i f\|_\infty  + \sum_{i,j}\|\partial_i\partial_j f\|_\infty$.
\item For a metric space $M$, $\cB(M)$ denotes the Borel $\sigma$-field on $M$.
\item For $n\ge 1$, $\cM(\RR^n)$ (resp. $\cM_b(\RR^n)$) denotes the space of measurable (resp. bounded measurable) functions $f:\RR^n\to\RR$.
\item We denote by $F$ the space $\cM(\RR)$. It will be equipped with the $\sigma$-field generated by $f\mapsto f(x)$ for all $x\in\RR$.
\item $\cP(\RR)$ denotes the space of probability measures on $(\RR,\cB(\RR))$. The space $\cP(\RR)$ is  equipped with the topology of narrow convergence. For $f\in\cM_b(\RR)$ and $\mu\in\cP(\RR)$, $\mu f$ or $\mu(f)$ denotes $\int_\RR f d\mu=\int_\RR f(x)\mu(dx)$.
\item A kernel is a measurable function $K$ from $\RR$ into $\cP(\RR)$. Denote by $E$ the space of all kernels on $\RR$. 
For $f\in \cM_b(\RR)$, $Kf\in \cM_b(\RR)$ is defined by $Kf(x)=\int_\RR f(y) K(x,dy)$. 
For $\mu\in\cP(\RR)$, $\mu K\in\cP(\RR)$ is defined by $(\mu K)f=\mu (Kf)$.  
If $K_1$ and $K_2$ are two kernels then $K_1K_2$ is the kernel defined by $(K_1K_2)f(x)=K_1(K_2f)(x)(=\int f(z)K_1(x,dy)K_2(y,dz))$. The space $E$ will be equipped with $\cE$ the $\sigma$-field generated by the mappings $K\mapsto \mu K$, for every $\mu\in\cP(\RR)$.
\item We denote by $\Delta$ (resp. $\Delta^{(n)}$ for $n\ge 1$) the Laplacian on $\RR$ (resp. on $\RR^n$), acting on twice differentiable functions $f$ on $\RR$ (resp. on $\RR^n$) and defined by $\Delta f=f^{''}$ (resp. $\Delta^{(n)}f=\sum_{i=1}^n \frac{\partial^2}{\partial x_i^2}f$).

\end{itemize}

\subsection{Definitions: Stochastic flows and $n$-point motions.}
\begin{definition} \label{defsfm} A measurable stochastic flow of mappings (SFM) $\p$ on $\RR$, defined on a probability space $(\Omega,\cA,\P)$, is a family $(\p_{s,t})_{s<t}$ such that
\begin{enumerate}
\item For all $s<t$, $\p_{s,t}$ is a measurable mapping from $(\Omega\times\RR,\cA\otimes\cB(\RR))$ to $(\RR,\cB(\RR))$;
\item For all $h\in\RR$, $s<t$, $\p_{s+h,t+h}$ is distributed like $\p_{s,t}$;
\item For all $s<t<u$ and all $x\in \RR$, a.s. $\p_{s,u}(x)=\p_{t,u}\circ \p_{s,t}(x)$, and $\p_{s,s}$ equals the identity;
\item For all $f\in C_0(\RR)$, and $s\le t$, we have
$$\lim_{(u,v)\to (s,t)}\sup_{x\in\RR}\E[(f\circ\p_{u,v}(x)-f\circ\p_{s,t}(x))^2]=0;$$
\item For all $f\in C_0(\RR)$, $x\in\RR$, $s<t$, we have
$$\lim_{y\to x}\E[(f\circ\p_{s,t}(y)-f\circ\p_{s,t}(x))^2]=0;$$
\item For all $s<t$, $f\in C_0(\RR)$, $\lim_{|x|\to\infty}\E[(f\circ\p_{s,t}(x))^2]=0$.
\end{enumerate}
\end{definition}

\begin{definition} \label{defsfk} A measurable stochastic flow of kernels (SFK) $K$ on $\RR$, defined on a probability space $(\Omega,\cA,\P)$, is a family $(K_{s,t})_{s<t}$ such that
\begin{enumerate}
\item For all $s<t$, $K_{s,t}$ is a measurable mapping from $(\Omega\times\RR,\cA\otimes\cB(\RR))$ to $(\cP(\RR),\cB(\cP(\RR)))$;
\item For all $h\in\RR$, $s<t$, $K_{s+h,t+h}$ is distributed like $K_{s,t}$;
\item For all $s<t<u$ and all $x\in \RR$, a.s. $K_{s,u}(x)=K_{s,t}K_{t,u}(x)$, and $K_{s,s}$ equals the identity;
\item For all $f\in C_0(\RR)$, and $s\le t$, we have
$$\lim_{(u,v)\to (s,t)}\sup_{x\in\RR}\E[(K_{u,v}f(x)-K_{s,t}f(x))^2]=0;$$
\item For all $f\in C_0(\RR)$, $x\in\RR$, $s<t$, we have
$$\lim_{y\to x}\E[(K_{s,t}f(y)-K_{s,t}f(x))^2]=0;$$
\item For all $s<t$, $f\in C_0(\RR)$, $\lim_{|x|\to\infty}\E[(K_{s,t}f(x))^2]=0$.
\end{enumerate}
\end{definition}

The law of a SFK (resp. of a SFM) is a probability measure on $(\Pi_{s<t} E,\otimes_{s\le t}\cE)$ (resp. on $(\Pi_{s<t} F,\otimes_{s\le t}\cF)$).
A SFK $K$ will be called a SFM when $K_{s,t}(x)=\delta_{\p_{s,t}(x)}$ for some SFM $\p$.

\begin{definition} \label{defnptsg} Let $(\P^{(n)}_t ,\;n\ge 1)$ be a family of Feller semigroups,
respectively defined on $\RR^n$ and acting on $C_0(\RR^n)$ as well
as on bounded continuous functions. We say that this family is
consistent as soon as
\begin{enumerate}
\item for all $n\ge 1$, all permutation $\sigma$ of $\{1,\dots,n\}$ and all $f\in C_0(\RR^n)$,
$$\P^{(n)}_t (f^\sigma)=(\P^{(n)}_tf)^\sigma$$ where $f^\sigma(x_1,\dots,x_n)=f(x_{\sigma_1},\dots,x_{\sigma_n})$;
\item for all $k \le n$, all $f\in C_0(\RR^k)$, all $y\in \RR^k$ and all $x\in \RR^n$ such that $x_i\le y_i$ for $i\le n$, we have
$$\P^{(n)}_tg(x)=\P^{(k)}_tf(y)$$
where $g(x_1,\dots,x_n)=f(y_1,\dots,y_k)$.
\end{enumerate}
We will denote by $\PP^{(n)}_x$ the law of the Markov process associated with
$\P^{(n)}_t$ and starting from $x\in\RR^n$. This Markov process will
be called the $n$-point motion of this family of semigroups.
\end{definition}



A general result (see Theorem 2.1 in \cite{LJR2004a}) states that there is a one to one correspondence between laws of SFK $K$ and consistent family of Feller semigroups $\P^n$, with the semigroup $\P^{(n)}$ defined by $\P^n_t=\E[K_{0,t}^{\otimes n}]$). It can be viewed as a generalization of De Finetti's Theorem (see \cite{Tsirelson04}).

\subsection{Definition: white noises}
Let $C:\RR\times\RR\to\RR$ be a covariance function on $\RR$, i.e. $C$ is symmetric and $\sum_{i,j}\lambda_i\lambda_jC(x_i,x_j)\ge 0$ for all finite family $(\lambda_i,x_i)\in\RR^2$. Assuming that the reproducing Hilbert space $H_C$ associated to $C$ is separable, there exists $(e_i)_{i\in I}$ (with $I$ at most countable) an orthonormal basis of $H_C$ such that $C(x,y)=\sum_{i\in I} e_i(x)e_i(y)$.

\begin{definition}\label{defwn}
A white noise of covariance $C$ is a centered Gaussian family of real random variables $(W_{s,t}(x),\; s<t,\; x\in\RR)$,  such that 

$$\E[W_{s,t}(x)W_{u,v}(y)]= |[s,t]\cap [u,v]| \times C(x,y).$$

A standard white noise is a centerd Gaussian family of real random variables $(W_{s,t},\; s<t)$ such that
$$\E[W_{s,t}(x)W_{u,v}(y)]= |[s,t]\cap [u,v]|.$$
\end{definition}

Starting with $(W^i)_{i\in I}$ independent standard white noises, one can define $W=(W_{s,t})_{s<t}$ a white noise of covariance $C$ by the formula
$$W_{s,t}(x)=\sum_{i\in I} W^i_{s,t} e_i(x),$$
which is well defined in $L^2$.

Although $W_{s,t}$ doesn't belong to $H_C$, one can recover $W^i_{s,t}$ out of $W$ by $W^i_{s,t}=\langle W_{s,t},e_i\rangle$. 
Indeed, for any given $i$, $e_i$ is the limit as $n\to\infty$ in $H_C$ of $e_i^n=\sum_k \lambda_k^n C_{x_k^n}$ where for all $n$, $(\lambda_k^n,x_k^n)$ is a finite family in $\RR^2$. Note that
$$\|e_i^n-e_i^m\|^2_{H_C}=\sum_{k,\ell}\lambda_k^n\lambda_\ell^m C(x_k^n,x_\ell^m).$$
Denote $W^{n,i}_{s,t}=\sum_k \lambda_k^n W_{s,t}(x_k^n)$. Then for $n$ and $m$,
$$\E[(W^{n,i}_{s,t}-W^{m,i}_{s,t})^2]= (t-s) \|e_i^n-e_i^m\|^2_{H_C}.$$
Thus one can define $W^i_{s,t}$ as the limit in $L^2$ of $W^{n,i}_{s,t}$. Then one easily checks that $W_{s,t}(x)=\sum_i W^i_{s,t}(x) e_i(x)$ a.s.

For $K$ a SFK (resp. $\p$ a SFM or $W$ a white noise), we denote for all $s\le t$ by $\cF^K_{s,t}$ (resp. $\cF^\p_{s,t}$ or $\cF^W_{s,t}$) the $\sigma$-field generated by $\{K_{u,v};~s\le u\le v\le t\}$ (resp. by $\{\p_{u,v};~s\le u\le v\le t\}$ or $\{W_{u,v};~s\le u\le v\le t\}$). A white noise $W$ is said to be a $\cF^K$ (resp. $\cF^\p$-white noise) if $\cF^W_{s,t}\subset \cF^K_{s,t}$ for all $s<t$ (resp. $\cF^W_{s,t}\subset \cF^\p_{s,t}$ for all $s<t$).
In all the following, all $\sigma$-fields will be completed by negligible events.

\subsection{Definition: the $(\frac{1}{2}\Delta,C)$-SDE}
Let $\p$ be a SFM and $W$ a $\cF^\p$-white noise of covariance $C$.
Then if for all $s<t$ and all $x\in\RR$, we have
\begin{equation} \label{edsCphi}
\p_{s,t}(x) = x + \sum_{i\in I} \int_s^t e_i \circ \p_{s,u}(x) W^i(du).
\end{equation}
then $(\p,W)$ is said to solve the SDE \eqref{edsCphi}.
Since this SDE is determined by the covariance $C$, we will more
simply say that $(\p,W)$ solves the $C$-SDE driven by $W$. Note that to find SFM's
solutions of the $C$-SDE for which the one-point motion is a Brownian
motion, we will need to assume that $C(x,x) = 1$ for all $x\in\RR$.
In all the following, we will be interested in constructing SFM and
SFK for which the one-point motion is a Brownian motion. Adding this
condition, the $C$-SDE will be called the
$(\frac{1}{2}\Delta,C)$-SDE, since $\frac{1}{2}\Delta f=\frac{1}{2}f"$ for
$f\in C^2(\RR)$ is the generator of the Brownian motion on $\RR$.

The notion of solution of this SDE can be extended to stochastic flows of kernels: let $K$ be a SFK and $W=\sum_i e_iW^i$, a $\cF^K$-white noise of covariance $C$, then $(K,W)$ is said to solve the $(\frac{1}{2}\Delta,C)$-SDE  driven by $W$ if 
\begin{equation}\label{edsC}
K_{s,t}f(x) = f(x) + \sum_{i\in I} \int_s^t K_{s,u} (e_i f')(x) W^i(du) +\frac{1}{2} \int_s^t K_{s,u}f''(x) du,
\end{equation}
for all $f\in C^2_0(\RR)$, $s<t$ and $x\in\RR$. 
Note that when $C$ is continuous, then identity \eqref{edsC} implies that $W$ is a $\cF^K$-white noise (see section 5 and Lemma 5.3 in \cite{LJR2004a}).
To find solutions of the $(\frac{1}{2}\Delta,C)$-SDE, we will need to assume that $C(x,x)\le 1$ for all $x\in\RR$. (This condition comes from the fact that for all $t\ge 0$,
$$\P_tf^2 = \E[K_{0,t}f^2(x)]\ge E[(K_{0,t}f)^2(x)] = \P^{(2)}_tf^{\otimes 2}f (x,x)$$
and that 
\begin{eqnarray*}
A f^2(x) - A^{(2)} f^{\otimes 2}(x,x) = (1-C(x,x)) (f'(x))^2
\end{eqnarray*} 
where $A=\frac{1}{2}\Delta$ is the generator of $\P_t$ and $A^{(2)}$ is the generator of $\P^{(2)}_t$.) Having $C$ stricly less than one means the flow is a mixture of stochastic transport and deterministic
heat flow.

Taking the expectation in \eqref{edsC}, we see that for
a solution $(K,W)$ of the $(\frac{1}{2}\Delta,C)$-SDE, the one-point motion of $K$ is a standard Brownian motion. Note that if $K$ is a SFK of the form $\delta_\p$, with $\p$ a SFM, then $(K,W)$ is a solution of the $(\frac{1}{2}\Delta,C)$-SDE if and only if $(\p,W)$ solves
\eqref{edsCphi}, and we must have $C(x,x)=1$ for all $x\in\RR$.

A solution $(K,W)$ of the $(\frac{1}{2}\Delta,C)$-SDE is called a Wiener solution  if for all $s<t$, $\cF^K_{s,t}=\cF^W_{s,t}$.
A SFM $\p$ will be called coalescing when for all $x,y$,
 $T=\inf\{t>0;\; \p_{0,t}(x)=\p_{0,t}(y)\}$ is finite a.s.
A SFK $K$ will be called diffusive when it is not a SFM.

\subsection{Martingale problems related to the $(\frac{1}{2}\Delta,C)$-SDE}
Let $(K,W)$ be a solution of the $(\frac{1}{2}\Delta,C)$-SDE. Denote
by $\P^{(n)}_t$, the semigroup associated with the $n$-point motion
of $K$, and by $\PP^{(n)}_x$ the law of the
$n$-point motion started from $x\in\RR^n$. 
\bprop Let $X^{(n)}$ be
distributed like $\PP^{(n)}_x$ with $x\in\RR^n$. Then it is a
solution of the martingale problem:
\begin{equation}
\label{MPC}
f(X^{(n)}_t)-\int_0^t A^{(n)}f(X^{(n)}_s)ds
\end{equation}
is a martingale for all $f\in C^2_0(\RR^n)$, where
$$A^{(n)}f(x)=\frac{1}{2}\Delta^{(n)} f(x) + \sum_{i<j} C(x_i,x_j)\frac{\partial^2}{\partial x_i\partial x_j}f(x).$$
\eprop \prf When $f$ is of the form $f_1\otimes \cdots \otimes f_n$,
with $f_1,\dots,f_n$ in $C_0^2(\RR)$, then it is easy to verify
\eqref{MPC} when $(K,W)$ solves \eqref{edsC}. This extends to the linear space spanned by such functions, and by density, to $C_0^2(\RR^n)$.
(A proof of this fact can be derived from the observation that, on the torus, Sobolev theorem shows that trigonometric polynomials are dense in $C^2$.) \qed

\begin{remark} If one can prove uniqueness for these martingale problems, then this implies there exists at most one solution to the $(\frac{1}{2}\Delta,C)$-SDE. (Indeed, using It\^o's formula, the expectation of any product of $K_{s,t}f(x)$ and $W_{s,t}(y)$'s can be expressed in terms of $\P^{(n)}_t$'s.)
\end{remark}

\subsection{Statement of the main Theorems}
\bthm \label{thm2}
Let $C_{\pm}(x,y)=1_{\{x>0\}}1_{\{y>0\}}+1_{\{x<0\}}1_{\{y<0\}}$. Then there is a unique solution $(K,W)$ of the $(\frac{1}{2}\Delta,C_\pm)$-SDE. Moreover,  this solution is a Wiener solution and $K$ is a coalescing SFM.
\ethm

Denote by $\p^{\pm}$ the coalescing flow defined in Theorem \ref{thm2}.
The white noise $W$ of covariance $C_{\pm}$ can be written in the form $W=W^+1_{\RR^+}+W^-1_{\RR^-}$, with $W^+$ and $W^-$ two independent white noises.
Then \eqref{edsC} is equivalent to
\begin{eqnarray}\label{edsCpm}
\p^{\pm}_{s,t}(x) &=& x + \int_s^t 1_{\{\p^{\pm}_{s,u}(x)>0\}}  W^+(du) + \int_s^t 1_{\{\p^{\pm}_{s,u}(x)<0\}}  W^-(du).
\end{eqnarray}
A consequence of this Theorem, with Proposition \ref{prop2b} below, is
\bthm \label{thm2b}
The SDE driven by  $B^+$ and $B^-$, two independent Brownian motions,
\begin{eqnarray}\label{edsCpmb}
dX_t = 1_{\{X_t>0\}} dB^+_t  + 1_{\{X_t<0\}} dB^-_t
\end{eqnarray}
has a unique solution. Moreover, this solution is a strong solution.
\ethm

\prf We recall that saying \eqref{edsCpmb} has a unique solution means that for all $x\in\RR$, there exists one and only one probability measure $\QQ_x$ on $C(\RR^+:\RR^3)$ such that under $\QQ_x(d\omega)$, the canonical process $(X,B^+,B^-)(\omega)=\omega$ satisfies \eqref{edsCpmb}, with $B^+$ and $B^-$ two independent Brownian motions, and $X$ a Brownian motion started at $x$. Proposition \ref{prop2b} states that \eqref{edsCpmb} has a unique solution. 
Since \eqref{edsCpm} holds, one can take $(X_t,B^+_t,B^-_t)=(\p^{\pm}_{0,t}(x),W^+_{0,t},W^-_{0,t})$ for this solution. 

A solution $(X,B^+,B^-)$ is a strong solution if $X$ is measurable with respect to the $\sigma$-field generated by $B^+$ and $B^-$, completed by the events of probability $0$. Since $\p^\pm$ is a Wiener solution, one can conclude. \qed

\bthm \label{thm1}
Let $C_+(x,y)=1_{\{x>0\}}1_{\{y>0\}}$.
\bdes \iti There is a unique solution $(K^+,W)$ solution of the $(\frac{1}{2}\Delta,C_+)$-SDE. 
\itii The flow $K^+$ is  diffusive and is a Wiener solution.
\itiii The flow $K^+$ can be obtained by filtering $\p^{\pm}$ with respect to the noise generated by $W^+$: 
$K^+_{s,t}=\E[\delta_{\p_{s,t}^\pm}|W^+]$.
\edes
\ethm

\bthm \label{thm4} There exists a unique coalescing SFM $\p^+$ such that its $n$-point motions coincide with the $n$-point motions of $K^+$ before hitting $\Delta_n=\{x\in\RR^n; \quad \exists i\neq j, \; x_i=x_j\}$. Moreover
\bdes
\iti $(\p^+,W^+)$ is not a solution of the $(\frac{1}{2}\Delta,C_\pm)$-SDE but it  satisfies
$$\int_s^t 1_{\{\p^{+}_{s,u}(x)>0\}} d\p^{+}_{s,u}(x) = \int_s^t 1_{\{\p^{+}_{s,u}(x)>0\}}  W^+(du).$$
\itii $K^+$ can be obtained by filtering $\p^+$ with respect to the noise generated by $W^+$: 
$K^+_{s,t}=\E[\delta_{\p_{s,t}^+}|W^+]$.
\edes
\ethm

\begin{remark} Following \cite{LJR2004a}, it should be possible to prove that he linear part of the noise generated by $\p^+$ is the noise generated by $W^+$.
\end{remark}

We refer to section \ref{filtering} for more precise definitions of noises, extension of noises, filtering by a subnoise, and linear part of a noise.



\section{General results.}
\subsection{Chaos decomposition of Wiener solutions.}
\bprop \label{21} Let $(K,W)$ be a Wiener solution of the $(\frac{1}{2}\Delta,C)$-SDE.
Then for all $s<t$, $f$ a bounded measurable function on $\RR$ and $x\in\RR$, a.s.,
\beq\label{chaosexp} K_{s,t}f(x)=\sum_{n\ge 0} J^n_{s,t}f(x) \eeq
with $J^n$ defined by $J^0_t=\P_t$, where $\P_t$ is the heat semigroup on $\RR$, and for $n\ge 0$,
\beq J^{n+1}_{s,t}f(x)=\sum_i\int_s^t J^n_{s,u}((\P_{t-u}f)'e_i)(x) W^i(du). \eeq
This implies that there exists at most one Wiener solution to the $(\frac{1}{2}\Delta,C)$-SDE.
\eprop
\begin{proof}
We essentially follow the proof of Theorem 3.2 in \cite{LJR2002}, with a minor correction at the end noticed by Bertrand Micaux during his PhD \cite{Micaux}.

Let us first remark that the stochastic integral $\sum_i \int_0^t K_s(e_i (\P_{t-s}f)')(x) W^i(ds)$ do converge in $L^2$ for all bounded measurable function $f$
since (using in the fourth inequality that $C(x,x)\le 1$ for all
$x\in\RR$) \beqarr
&&\hspace{-50pt}\sum_i\E\left[\left(\int_s^t K_{s,u}((\P_{t-u}f)'e_i)(x)W^i(du)\right)^2\right]\\
&\le& \sum_i \int_s^t\E\left[\left(K_{s,u}((\P_{t-u}f)'e_i)(x)\right)^2\right]du\\
&\le& \sum_i \int_s^t\E\left[K_{s,u}\big((\P_{t-u}f)'e_i\big)^2(x)\right]du\\
&\le& \int_s^t\P_{u-s}\left(\big((\P_{t-u}f)'\big)^2\sum_i e_i^2\right)(x) du\\
&\le& \int_s^t\P_{u-s}\big((\P_{t-u}f)'\big)^2(x) du
\eeqarr
which is finite (since, using that $\frac{d}{du} \big(\P_{u-s}(\P_{t-u}f)^2\big) = \P_{u-s}\big((\P_{t-u}f)'\big)^2$, it is equal to $\P_{t-s}f^2(x)-(\P_{t-s}f)^2(x)$).

Take $f$ a $C^3$ function with compact support.
For $t>0$, denote $K_{0,t}$ simply $K_t$.
Then for $t>0$, $n\ge 1$ and $x\in\RR$
\beqarr
K_tf(x)-\P_tf(x)
&=& \sum_{k=0}^{n-1}  \left(K_{t(k+1)/n} (\P_{t(1-(k+1)/n)}f) - K_{tk/n} (\P_{t(1-k/n)}f)\right) \\
&=& \sum_{k=0}^{n-1}  (K_{t(k+1)/n}-K_{tk/n}) (\P_{t(1-(k+1)/n)}f) \\
&& +  \sum_{k=0}^{n-1} K_{tk/n} (\P_{t(1-(k+1)/n)}f-\P_{t(1-k/n)}f) \\
&=&  \sum_{k=0}^{n-1} \sum_i\int_{tk/n}^{t(k+1)/n}  K_s\left(e_i(\P_{t(1-(k+1)/n)}f)'\right)(x) W^i(ds)\\
&& +  \sum_{k=0}^{n-1} \int_{tk/n}^{t(k+1)/n} K_s\left(\frac{1}{2}(\P_{t(1-(k+1)/n)}f)''\right)(x) ds \\
&& -  \sum_{k=0}^{n-1} K_{tk/n} (\P_{t(1-k/n)}f-\P_{t(1-(k+1)/n)}f) 
\eeqarr
since $(K,W)$ solves the $(\frac{1}{2}\Delta,C)$-SDE.
This last expression implies that for $t>0$, $n\ge 1$ and $x\in\RR$,
$$K_tf(x)-\P_tf(x)-\sum_i \int_0^t K_s(e_i (\P_{t-s}f)')(x) W^i(ds) = \sum_{k=1}^3 B_k(n)$$
with
\beqarr
B_1(n) &=& \sum_{k=0}^{n-1} \sum_i\int_{tk/n}^{t(k+1)/n}  K_s\left(e_i(\P_{t(1-(k+1)/n)}f-\P_{t-s}f)'\right)(x) W^i(ds),\\
B_2(n) &=& -\sum_{k=0}^{n-1} K_{tk/n} \left(\P_{t(1-k/n)}f-\P_{t(1-(k+1)/n)}f - \frac{t}{2n} (\P_{t(1-(k+1)/n)}f)''\right)(x),\\
B_3(n) &=& \sum_{k=0}^{n-1} \int_{tk/n}^{t(k+1)/n} (K_s-K_{tk/n})\left(\frac{1}{2}(\P_{t(1-(k+1)/n)}f)''\right)(x) ds.
\eeqarr
The terms in the expression of $B_1(n)$ being orthogonal,
\beqarr
\E[(B_1(n))^2]
&=& \sum_{k=0}^{n-1} \sum_i\int_{tk/n}^{t(k+1)/n} \E\left[\left(K_s\big(e_i(\P_{t(1-(k+1)/n)}f-\P_{t-s}f)'\big)\right)^2(x)\right] ds\\
&\le& \sum_{k=0}^{n-1} \sum_i\int_{tk/n}^{t(k+1)/n} \P_s\left(e_i(\P_{t(1-(k+1)/n)}f-\P_{t-s}f)'\right)^2(x) ds\\
&\le& \sum_{k=0}^{n-1} \int_{tk/n}^{t(k+1)/n} \P_s\left((\P_{t(1-(k+1)/n)}f-\P_{t-s}f)'\right)^2(x) ds
\eeqarr
where we have used Jensen inequality in the second inequality and the fact that $\sum_i e_i^2(x)=C(x,x)\le 1$ in the last inequality.
This last term is less that $n\int_0^{t/n}\|\P_sf'-f'\|^2_\infty ds=O(\|f''\|^2_\infty t^{2}/ n)$. 

By using triangular inequality, $\E[(B_2(n))^2]^{1/2}$ is less than
\beqarr
&&\hskip-20pt \sum_{k=0}^{n-1} \E\left[K_{tk/n} \left(\P_{t(1-k/n)}f-\P_{t(1-(k+1)/n)}f - \frac{t}{2n} (\P_{t(1-(k+1)/n)}f)''\right)^2(x)\right]^{1/2}\\
&\le& \sum_{k=0}^{n-1} \left[\P_{tk/n} \left(\P_{t(1-k/n)}f-\P_{t(1-(k+1)/n)}f - \frac{t}{2n} (\P_{t(1-(k+1)/n)}f)''\right)^2(x)\right]^{1/2}.\\ 
\eeqarr
Using moreover that $\|\P_{t/n}f-f-t/(2n)f''\|_\infty=O((t/n)^{3/2}\|f'''\|_\infty)$, we get that $\E[(B_2(n))^2]^{1/2}=O(t^{3/2}n^{-1/2}\|f'''\|_\infty)$.

By using again triangular inequality, $\E[(B_3(n))^2]^{1/2}$ is less than
\beqarr
&&\sum_{k=0}^{n-1} \E\left[\left(\int_{tk/n}^{t(k+1)/n} (K_s-K_{tk/n})\left(\frac{1}{2}(\P_{t(1-(k+1)/n)}f)''\right)(x) ds\right)^2\right]^{1/2}\\
&\le& \sum_{k=0}^{n-1} \left(\frac{t}{n}\int_{tk/n}^{t(k+1)/n} \E\left[\left((K_s-K_{tk/n})\left(\frac{1}{2}(\P_{t(1-(k+1)/n)}f)''\right)\right)^2(x)\right] ds\right)^{1/2}.
\eeqarr
For $f$ a Lipschitz function and $0\le s<t$, we have (with $(X,Y)$ the two point motion of $K$ of law  $\PP^{(2)}_{(x,x)}$ and $\E^{(2)}_{(x,x)}$ the expectation with respect to $\PP^{(2)}_{(x,x)}$)
\beqarr
\E[(K_tf-K_sf)^2(x)]
&=& \E^{(2)}_{(x,x)}[(f(X_{t})-f(X_s))(f(Y_{t})-f(Y_s))]\\
&\le& (t-s) (\hbox{Lip}(f))^2
\eeqarr
with $\hbox{Lip}(f)$ the Lipschitz constant of $f$.
From this estimate, and since $\hbox{Lip}(\P_{t(1-(k+1)/n)}f'')\le\hbox{Lip}(f'')$, we deduce that $\E[(B_3(n))^2]^{1/2}=O(t^{3/2}\hbox{Lip}(f'')n^{-1/2})$.
The estimates we gave for $\E[(B_i(n))^2]$, for $i\in\{1,2,3\}$, implies that for $f$ a $C^3$ function with compact support,
$$K_tf(x) = \P_tf(x) + \sum_i \int_0^t K_s(e_i (\P_{t-s}f)')(x) W^i(ds).$$ 

This implies that \beq\label{relchaos}
K_{s,t}f(x)=\P_{t-s}f(x)+\sum_i\int_s^t K_{s,u}((\P_{t-u}f)')(x)
W^i(du).\eeq 

Iterating relation \eqref{relchaos}, we conclude from the orthogonality of Wiener chaoses that \eqref{chaosexp} holds for $f$ a $C^3$ function with compact support. This extends to all bounded measurable functions. 
 \end{proof}

\subsection{Filtering a SFK by a subnoise.} \label{filtering}
We follow here section 3 in \cite{LJR2004a}.
\begin{definition}\label{def:noise} A noise consists of a separable probability space $(\Omega,\cA,\P)$,
a one-parameter group $(T_h)_{h\in\RR}$ of $\P$-preserving $L^2$-continuous transformations
of $\Omega$ and a family ${\cF_{s,t},\;\le s \le t \le \infty}$ of sub-$\sigma$-fields of $\cA$ such that:
\bdes
\ita $T_h$ sends $\cF_{s,t}$ onto $\cF_{s+h,t+h}$ for all $h \in \RR$ and $s \le t$ ,
\itb $\cF_{s,t}$ and $\cF_{t,u}$ are independent for all $s \le t \le u$,
\itc $\cF_{s,t} \wedge \cF_{t,u} = \cF_{s,u}$ for all $s \le t \le u$.
\edes
Moreover, we will assume that, for all $s \le t$ , $\cF_{s,t}$ contains all $\P$-negligible sets
of $\cF_{-\infty,\infty}$, denoted $\cF$ .
\end{definition}

A subnoise $\bar{N}$ of $N$ is a noise $(\Omega,\cA,\P,(T_h),\bar{\cF}_{s,t})$, with $\bar{\cF}_{s,t}\subset \cF_{s,t}$. Note that a subnoise is characterized by  a $\sigma$-field invariant by $T_h$ for all $h\in\RR$, $\bar{\cF}=\bar{\cF}_{-\infty,\infty}$. In this case, we will say that $\bar{N}$ is the noise generated by $\bar{\cF}$.

A linear representation of $N$ is a family of real random variables
$X = (X_{s,t};\; s \le t)$ such that:
\begin{description}
\ita  $X_{s,t}\circ T_h = X_{s+h,t+h}$ for all $s \le t$ and $h\in\RR$,
\itb $X_{s,t}$ is $\cF_{s,t}$-measurable for all $s \le t$ ,
\itc $X_{r,s} +X_{s,t} = X_{r,t}$ a.s., for all $r \le s \le t$.
\end{description}
Define $\cF^{lin}$ be the $\sigma$-field generated by the random variables $X_{s,t}$ where $X$ is a linear representation of $N$ and $s\le t$. Define $N^{lin}$ to be the subnoise of $N$ generated by $\cF^{lin}$. This noise is called the linear part of the noise $N$.

Let $\P^0$ be the law of a SFK, it is a law on $(\Omega^0,\cA^0)=(\Pi_{s<t}E,\otimes_{s<t}\cE)$, and let $K$ be the canonical SFK of law $\P^0$. For $h\in\RR$, define $T^0_h:\Omega\to\Omega$ by $T^0_h(\omega^0)_{s,t}=\omega^0_{s+h,t+h}$. For $s<t$, $\cF^0_{s,t}$ is the $\sigma$-field $\cF^K_{s,t}$ completed by $\P^0$-negligible sets of $\cA^0$. This defines a noise $N^0$ called the noise of the SFK $K$.

Let $\bar{N}$ be a subnoise of $N^0$.
In section 3.2 in \cite{LJR2004a}, a SFK $\bar{K}$ is defined as the filtering of $K$ with respect to $\bar{N}$: for $s<t$, $f\in C_0(\RR)$ and $x\in \RR$, $\bar{K}$ is such that
$$\bar{K}_{s,t}f(x)=\E[K_{s,t}f(x)|\bar{\cF}_{s,t}]=\E[K_{s,t}f(x)|\bar{\cF}].$$

Suppose now $\PP_0$ is the law of a SFK that solves the $(\frac{1}{2}\Delta,C)$-SDE. Writing $C$ in the form $C(x,y)=\sum_{i\in I} e_i(x)e_i(y)$, there is a $\cF^K$-white noise $W=\sum_{i\in I}W^ie_i$ such that:
\begin{equation}\label{sdeI}K_{s,t}f(x) = f(x) + \sum_{i\in I} \int_s^t K_{s,u} (e_i f')(x) W^i(du) +\frac{1}{2} \int_s^t K_{s,u}f''(x) du.\end{equation}
(i.e. the $(\frac{1}{2}\Delta,C)$-SDE is transported on the canonical space).
\bprop \label{prop:edsfilt} Let $J\subset I$ and $\bar{N}$ the noise generated by $\{W^j; j\in J\}$, and $\bar{W}=\sum_{j\in J}W^j e_j$. Let $\bar{K}$ be the SFK obtained by filtering $K$ with respect to $\bar{N}$. Then \eqref{sdeI} is satisfied with $K$ replaced by $\bar{K}$ and $W$ by $\bar{W}$ (or $I$ by $J$). If $\sigma(\bar{W})\subset \sigma(\bar{K})$, then 
$\bar{K}$ solves the $(\frac{1}{2}\Delta,\bar{C})$-SDE  with $\bar{C}(x,y)=\sum_{j\in J} e_j(x)e_j(y)$, and driven by $\bar{W}$.
\eprop
\prf This proposition reduces to prove that
$$\bar{K}_{s,t}f(x) = f(x) + \sum_{j\in J} \int_s^t \bar{K}_{s,u} (e_j f')(x) W^j(du) +\frac{1}{2} \int_s^t \bar{K}_{s,u}f''(x) du $$
which easily follows by taking the conditional expectation with respect to $\sigma(W^j;\;j\in J)$ in equation \eqref{sdeI}. \qed

\section{$C_{\pm}(x,y)=1_{\{x<0\}}1_{\{y<0\}} + 1_{\{x>0\}}1_{\{y>0\}}$.}

\subsection{The SDE $dX_t=1_{\{X_t>0\}}dB^{+}_t + 1_{\{X_t<0\}}dB^{-}_t$.}
Let $B^1$ and $B^2$ be two independent Brownian motions.
For $x\in\RR$, define $X^x$, $B^{x,+}$ and $B^{x,-}$ by
\begin{eqnarray}
X^x_t &=& x +B^1_t\\
B^{x,-}_t &=& \int_0^t 1_{\{X^x_s<0\}}dX^x_s+\int_0^t 1_{\{X^x_s>0\}}dB^2_s.\\
B^{x,+}_t &=& \int_0^t 1_{\{X^x_s<0\}}dB^2_s+\int_0^t 1_{\{X^x_s>0\}}dX^x_s.
\end{eqnarray}
Then $X^x$, $B^{x,-}$ and $B^{x,+}$ are Brownian motions respectively started at $x$, $0$ and $0$. Moreover $B^{x,-}$ and $B^{x,+}$ are independent, and we have
\beq X_t^x=x+\int_0^t 1_{\{B^x_s>0\}}dB^{x,+}_s + \int_0^t 1_{\{B^x_s<0\}}dB^{x,-}_s. \eeq
Denote by $\QQ_x$ the law of $(X^x,B^{x,-},B^{x,+})$.


\bprop \label{prop2b}  Let $x\in \RR$, $X$, $B^+$ and $B^-$ be real random processes such that $B^+$ and $B^-$ are independent Brownian motions. Then if
\beq X_t=x+\int_0^t 1_{\{X_s>0\}}dB^{+}_s + \int_0^t 1_{\{X_s<0\}}dB^{-}_s, \eeq
the process $(X,B^+,B^-)$ has for law $\QQ_x$.
\eprop
\begin{proof} Let $$B_t=\int_0^t 1_{\{X_s<0\}}dB^{+}_s + 1_{\{X_s>0\}}dB^{-}_s.$$
Observe that $X-x$ and $B$ and two independent Brownian motions. Moreover
\begin{eqnarray}
B^{+}_t &=& \int_0^t 1_{\{X_s<0\}}dB_s+\int_0^t 1_{\{X_s>0\}}dX_s\\
B^{-}_t &=& \int_0^t 1_{\{X_s<0\}}dX_s + \int_0^t 1_{\{X_s>0\}}dB_s.
\end{eqnarray}
This implies the proposition ($(X,B^+,B^-)$ is distributed like $(X^x,B^{x,+},B^{x,-})$). \end{proof}


\subsection{Construction of the $n$-point motions up to $T^{(n)}$.} \label{sec:npt}
Let $x\not\in\Delta_n=\{x\in\RR^n;\; \exists i\neq j,\; x_i=x_j\}$.
For convenience, we will assume that $x_1<x_2<\cdots<x_n$ and set $i$ the integer such that $x_i\le 0 <x_{i+1}$, with $i=1$ when $x_1>0$ and $i=n$ when $x_n\le 0$. Let $B^1$ and $B^2$ be two independent Brownian motions.
In the following construction, $X_i$  follows $B^1$. Out of $X_i$ and $B^2$, we construct $B^+$ and $B^-$ two independent Brownian motions, and for $j>i$ (resp. for $j<i$), $X_j$ follows $B^+$ (resp. $B^-$), this until the first time $\tau$ when $X_{i+1}$ or $X_{i-1}$ hits $0$. After time $\tau$, we follow the same procedure by replacing $i$ by $i-1$ when $X_{i-1}(\tau)=0$ or by $i+1$ when $X_{i+1}(\tau)=0$. More precisely,
define for $t>0$, the processes
\beqarr {X}^0_i(t) &=& x_i + B^1_t\\
{B}^{0,-}_t &=&\int_0^t 1_{\{{X}^0_i(s)<0\}}dB^1_s + \int_0^t 1_{\{{X}^0_i(s)>0\}}dB^2_s  \\
{B}^{0,+}_t &=& \int_0^t 1_{\{{X}^0_i(s)<0\}}dB^2_s+ \int_0^t 1_{\{{X}^0_i(s)>0\}}dB^1_s \\
{X}^0_j(t) &=& x_j+ {B}^{0,-}_t \qquad\hbox{ for } j\le i-1 \\
{X}^0_j(t) &=& x_j+{B}^{0,+}_t \qquad\hbox{ for } j\ge i+1.
\eeqarr
Set $$\tau_1=\inf\{t>0: {X}^0_{i-1}(t)=0 \hbox{ or } {X}^0_{i+1}(t)=0\}$$ and set for $t\leq \tau_1$, $({X},{B}^-,{B}^+)(t)=({X}^0,{B}^{0,-},{B}^{0,+})(t)$.
Set $i_1=i+1$ if $X_{i+1}(\tau_1)=0$ and $i_1=i-1$ if $X_{i-1}(\tau_1)=0$

Then ${X}_1(\tau_1)<{X}_2(\tau_1)<\cdots<X_{i_1}(\tau_1)=0<\cdots
<X_n(\tau_1)$. Assume now that $(\tau_k)_{k\le \ell}$ and
$({X},{B}^-,{B}^+)(t)$  have been defined for $t\le\tau_\ell$ such
that a.s.
\begin{itemize}
\item $(\tau_k)_{1\le k\le\ell}$ is an increasing sequence of stopping times with respect to the filtration associated to ${X}$;
\item ${X}_1(\tau_k)<\cdots<{X}_n(\tau_k)$ for $1\le k\le \ell$;
\item for all $k$, there exists an integer $i_k$ such that ${X}_{i_k}(\tau_k)=0$.
\end{itemize}
We then define $(X_t,{B}^-_t,{B}^+_t)_{t\in ]\tau_\ell, \tau_{\ell+1}]}$ as is defined $({X}_t,{B}^-_t,{B}^+_t)_{t\le \tau_{1}}$ by replacing $i$ by $i_\ell$, $x$ by ${X}_{\tau_\ell}$ and $(B^1_\cdot,B^2_\cdot)$ by $(B^1_{\tau_\ell+\cdot}-B^1_{\tau_\ell},B^2_{\tau_\ell+\cdot}-B^2_{\tau_\ell})$.
Let ${T}=\lim_{\ell\to\infty}\tau_\ell$. Note that $T=\inf\{t\ge 0;\; X_t\in\Delta_n\}$ (with the convention $\inf\emptyset=+\infty$).
Denote by $\PP^{(n),0}_x$ the law of $(X_t)_{t\leq T}$.

\blem \label{lem:MPunicity} Let $X^{(n)}$ be a solution to the martingale problem \eqref{MPC}, with $X^{(n)}_0=x$. Let $T^{(n)}=\inf\{t;\; X^{(n)}_t\in \Delta_n)\}$. Then  $(X^{(n)}_t)_{t\le T^{(n)}}$ is distributed like $\PP^{(n),0}_x$. \elem
\prf Let $x=(x_1,\dots,x_n)\in\RR^n$ be such that $x_1<\cdots<x_n$.
Again, let $i$  be such that $x_i\le 0 <x_{i+1}$, with $i=1$ when $x_1>0$ and $i=n$ when $x_n\le 0$. Let $X^{(n)}$ be a solution of the martingale problem.
This implies that for all $j$, $X^{(n)}_j$ is a Brownian motion and for all $j$ and $k$,
\begin{equation}
\label{crochet}\langle X^{(n)}_j,X^{(n)}_k\rangle_t = \int_0^t C_{\pm}(X^{(n)}_j(s),X^{(n)}_k(s)) ds.\end{equation}

Let $B^0$  be a Brownian motion, independent of $X^{(n)}$.
Set $$\tau_1=\inf\{t>0: \exists j\neq i,\quad X^{(n)}_{j}(t)=0\}.$$
Define for $t\le \tau_1$,
\begin{eqnarray*}
B^+_t &=& X^{(n)}_{i+1}(t)-x_{i+1}, \qquad \hbox{when } i\le n-1\\
B^+_t&=& \int_0^t 1_{\{{X}^{(n)}_i(s)>0\}}d{X}^{(n)}_i(s)  + \int_0^t 1_{\{{X}^{(n)}_i(s)<0\}}dB^0_s, \qquad \hbox{when } i=n.
\end{eqnarray*}
Define also for $t\le \tau_1$
\begin{eqnarray*}
B^-_t &=& X^{(n)}_{i-1}(t)-x_{i-1}, \qquad \hbox{when } i\ge 2\\
B^-_t&=& \int_0^t 1_{\{{X}^{(n)}_i(s)<0\}}d{X}^{(n)}_i(s)  + \int_0^t 1_{\{{X}^{(n)}_i(s)>0\}}dB^0(s), \qquad \hbox{when } i=1.
\end{eqnarray*}
Define for $t\le \tau_1$,
$$B^1_t=X^{(n)}_i(t)-x_i$$
and
$$B^2_t=\int_0^t 1_{\{{X}^{(n)}_i(s)>0\}}dB^-_t + \int_0^t 1_{\{{X}^{(n)}_i(s)<0\}}dB^+_t.$$

Note that for $t\le \tau_1$, $\langle B^1,B^2\rangle_t=0$ and that
\beqarr
{X}^{(n)}_j(t)&=& x_j+B^-_t \qquad\hbox{ for } \quad j<i\\
{X}^{(n)}_j(t) &=& x_j+{B}^{+}_t \qquad \hbox{ for } \quad j>i.
\eeqarr

Assume now that $(\tau_k)_{k\le \ell}$ and $(B^1_t,B^2_t,{B}^+_t,B^-_t)_{t\le\tau_\ell}$ have been defined such that a.s.
\begin{itemize}
\item $(\tau_k)_{1\le k\le\ell}$ is an increasing sequence of stopping times with respect to the filtration associated to ${X^{(n)}}$;
\item $X^{(n)}_{\tau_k}\not\in \Delta_n$ for $1\le k\le \ell$;
\item for all $1\le k\le \ell$, there exists an integer $i_k$ such that $X^{(n)}_{i_k}(\tau_k)=0$.
\end{itemize}
We then define $(B^1_t,B^2_t,B^+_t,B^-_t)_{\tau_\ell<t\le \tau_{\ell+1}}$ as is defined $(B^1_t,B^2_t,B^+_t,B^-_t)_{0<t\le \tau_{1}}$ by replacing $i$ by $i_\ell$, $x$ by ${X^{(n)}}(\tau_\ell)$.
Note that ${T^{(n)}}=\lim_{\ell\to\infty}\tau_\ell$. Define for $t\ge T^{(n)}$, $B^1_t=B^1_{T^{(n)}}+X^{(n)}_1(t)-X^{(n)}_1(T^{(n)})$ and $B^2_t=B^2_{T^{(n)}}+B^0_t-B^0_{T^{(n)}}$. Then $B^1$ and $B^2$ are two independent Brownian motions (it is also the case for $B^+$ and $B^-$).
We finally remark that $(X^{(n)}_t)_{t\le T^{(n)}}$ can be defined out of $B^1$ and $B^2$ in the same way $\PP^{(n),0}_x$ is defined. This proves the lemma.
\qed

\subsection{Brownian motions with oblique reflection.}
Let $X$ and $B$ be two independent Brownian motions, with $X$ started at $x\in\RR$. Let $Y$ be a Brownian motion started at $y$ defined by
$$Y_t=y+\int_0^t1_{\{X_s<0\}}dB_s +\int_0^t1_{\{X_s>0\}}dX_s.$$
Let $A_t= \int_0^t 1_{\{X_s<0\}} ds$ and $\kappa_t=\inf\{ s>0;\; A_s>t\}$. 
Then define $(X^r,Y^r)_t=(X,Y)_{\kappa_t}$, it is a continuous process taking its values in $\{x\le 0\}$. Denote by $L_t$ the local time at $0$ of $X_t$.
Then $X^r$ is a Brownian motion in $\RR^-$ reflected at $0$ and that if $L^r_t=\frac{1}{2}L_{\kappa_t}$, $B^1_t=X_t^r-L_t^r-X_0^r$ is a Brownian motion.
\blem \label{lemoblique}  Set $B^2_t=\int_0^{\tau_t} 1_{\{X_s<0\}}dB_s$. Then $B^2$ is a Brownian motion and we have
$$Y^r_t=B^2_t-L^r_t+Y^r_0.$$
Thus $(X^r,Y^r)$ is a Brownian motion in $\{x\le 0\}$ reflected at $\{x=0\}$ with oblique reflection with angle of reflection equal to $\pi/4$, i.e.
$$(X^r_t,Y^r_t)=(X^r_0,Y^r_0)+ (B^1_t,B^2_t)-(L^r_t,L^r_t).$$
(note that $(X^r_0,Y^r_0)=(x,y)$ if $x\le 0$ and $(X^r_0,Y^r_0)=(0,y-x)$ if $x> 0$).
\elem
\prf Note that  $B^1_t=\int_0^{\tau_t} 1_{\{X_s<0\}}dX_s$. Since $\langle B^1\rangle_t=\langle B^2\rangle_t=\int_0^{\tau_t}1_{\{X_s<0\}} ds$ and since $\langle B^1,B^2\rangle_t=0$, $B^1$ and $B^2$ are two independent Brownian motions. Let $\eps$ be a small positive parameter and define the sequences of stopping times $\sigma_k^\eps$ and $\tau_k^\eps$ such that $\sigma_0^\eps=0$ and for $k\ge 0$
\beqarr \tau_k^\eps &=& \inf\{t\ge \sigma_k^\eps; X_t=0\},\\
\sigma_{k+1}^\eps &=& \inf\{t\ge \tau_k^\eps; X_t=\eps\}.
\eeqarr

Note that $X^r_0=x$ if $x\le0$ and $X^r_0=0$ if $x>0$, and for $t>0$
$$X^{r}_t =  \sum_{k\ge 1} (X_{\tau^\eps_k\wedge \kappa_t}-X_{\tau^\eps_{k-1}\wedge \kappa_t}) +X_{\tau_0^\eps\wedge\kappa_t},$$
with $X_{\tau_0^\eps\wedge\kappa_t}=X_t$ if $x\le 0$ and $\kappa_t<\tau_0^\eps$ and $X_{\tau_0^\eps\wedge\kappa_t}=0$ if  $\kappa_t>\tau_0^\eps$ (which holds if $x> 0$).
We have that $X^{r}_t=-L^{\eps,r}_t+B^{\eps,1}_t+X^r_0$ with
\beqarr L^{\eps,r}_t &=& -\sum_{k\ge 1} (X_{\tau^\eps_k\wedge \kappa_t}-X_{\sigma^\eps_{k}\wedge \kappa_t})\\
B^{\eps,1}_t &=& \sum_{k\ge 1} (X_{\sigma^\eps_k\wedge \kappa_t}-X_{\tau^\eps_{k-1}\wedge \kappa_t})+(X_{\tau_0^\eps\wedge\kappa_t}-X^r_0).
\eeqarr
Then $L^{\eps,r}_t$ and $B^{\eps,1}_t$ both converges in probability respectively towards $L^r_t$ and $B^1_t$ (see Theorem 2.23 in chapter 6 of \cite{KS}).

Note now that $Y^r_0=y$ is $x\le 0$ and $Y^r_0=y-x$ if $x>0$, and for $t>0$
$$Y^{r}_t =  \sum_{k\ge 1} (Y_{\tau^\eps_k\wedge \kappa_t}-Y_{\tau^\eps_{k-1}\wedge \kappa_t}) + Y_{\tau_0^\eps\wedge\kappa_t},$$
with $Y_{\tau_0^\eps\wedge\kappa_t}=Y_t$ if $x\le 0$ and $\kappa_t<\tau_0^\eps$ and $Y_{\tau_0^\eps\wedge\kappa_t}=y-x$ if  $\kappa_t>\tau_0^\eps$ (which holds if $x> 0$).
Since when $X_t$ and $Y_t$ are positive, $X_t-Y_t$ remains constant,
$Y^{r}_t=-L^{\eps,r}_t+B^{\eps,2}_t+Y^r_0$, with
$$B^{\eps,2}_t = \sum_{k\ge 1} (Y_{\sigma^\eps_k\wedge \kappa_t}-Y_{\tau^\eps_{k-1}\wedge \kappa_t}) + (Y_{\tau_0^\eps\wedge\kappa_t}-Y^r_0).$$
It remains to observe that $B^{\eps,2}_t$ converges in probability towards $B^2_t$. \qed

\medskip
Denote now by $(X,Y)$ a solution of the martingale problem \eqref{MPC}, for $n=2$. Denote by $T=\inf\{t\ge 0;\; X_s=Y_s\}$. Let $D^+=\{(x,y)\in\RR^2;\;x\le 0 \hbox{ and } y\ge 0\}$, $D^-=\{(x,y)\in\RR^2;\;x\ge 0 \hbox{ and } y\le 0\}$ and $D=D^+\cup D^-$.
Let $A^{\pm}_t=\int_0^{t\wedge T} 1_{\{(X_s,Y_s)\in D^{\pm}\}}ds$ and $A_t=A^+_t+A^-_t$.
Let $\kappa^{\pm}_t=\inf\{s;\;A^{\pm}_s>t\}$ and $\kappa_t=\inf\{s;\;A_s>t\}$. Set $(X^{r,\pm}_t,Y^{r,\pm}_t)=(X,Y)(\kappa^{\pm}_t)$ and $(X^{r}_t,Y^{r}_t)=(X,Y)(\kappa_t)$.
 Note that if $X_0>Y_0$ (resp. if $X_0<Y_0$) we have $(X^{r},Y^{r})=(X^{r,+},Y^{r,+})$ (resp. $(X^{r},Y^{r})=(X^{r,-},Y^{r,-})$).
Denote by $L_t(X)$ and by $L_t(Y)$ the local times at $0$ of $X$ and of $Y$. And denote $\frac{1}{2}L_{\kappa_t}(x)$ and $\frac{1}{2}L_{\kappa_t}(Y)$ by $L^1_t$ and by $L^2_t$.

The proof of Lemma \ref{lemoblique} can be adapted to prove that
\blem
The process $(X^{r},Y^{r})$ is a Brownian motion in $D$, with oblique reflection at the boundary of angle of reflection equal to $\pi/4$, and stopped when it hits $(0,0)$. More precisely, there exists $B^1$ and $B^2$ two independent Brownian motions such that when $X_0<Y_0$,
\beqarr
X^r_t &=& X^r_0 + B^1_t - L^1_t + L^2_t\\
Y^r_t &=& Y^r_0 + B^2_t - L^1_t + L^2_t
\eeqarr
and when $X_0>Y_0$,
\beqarr
X^r_t &=& X^r_0 + B^1_t + L^1_t -L^2_t\\
Y^r_t &=& Y^r_0 + B^2_t + L^1_t - L^2_t
\eeqarr
with $(X^r_0,Y^r_0)=(X_0,Y_0)$ if $(X_0,Y_0)\in D$,  $(X^r_0,Y^r_0)=(0,Y_0-X_0)$ if $Y_0<X_0<0$ or if $0<X_0<Y_0$, and  $(X^r_0,Y^r_0)=(X_0-Y_0,0)$ if $0<Y_0<X_0$ or if $X_0<Y_0<0$.
\elem

This process is a special case of the class of processes studied in \cite{VW85,W85}: the process $(X^r,Y^r)$ is a Brownian motion in
the wedge $D$, with angle $\xi=\pi/2$ and angles of reflection
$\theta_1=\theta_2=\pi/4$. An important parameter in the study of
these processes is $\alpha=(\theta_1+\theta_2)/\xi=1$.

\medskip
The fact that the reflected Brownian motion is a semimartingale on $[0,\tau_0]$, where $\tau_0$ is the first time to hit the origin, follows from Theorem 1 of \cite{W85}. An alternative proof, in a more general multi-dimensional and state-dependent coefficient setting,
is given in  Theorem 1.4 (property 3) of \cite{Ramanan06}.

The non-semimartingale property of the Markovian extension spending
no time at zero follows from Theorem 5 (with $\alpha = 1$) of
\cite{W85}. Two alternative proofs, the first of which also
generalizes to more general diffusions and higher dimensions, are
given in Theorem 3.1 of \cite{Ramanan10} and Proposition 4.13 of
\cite{Ramanan09}. These result are essentially related to the study
of the local times which is done in the next section. We choose not
to derive the results of this study from these references in order
to keep their proof self contained and reasonably short.

\subsection{Coalescing $n$-point motions.}
Without loss of generality, assume that $(X^r_0,Y^r_0)\in D^+$. Let $U_t=\frac{Y^r_t+X^r_t}{\sqrt{2}}$ and $V_t=\frac{Y^r_t-X^r_t}{\sqrt{2}}$. Then $V$ is a Brownian motion stopped when it hits $0$. This implies in particular that $T_0=\inf\{t\ge 0; \; X^r_t=Y^r_t=0\}$ is finite a.s. Denote by $L^r:=L^1+L^2$ the local time at $0$ of $(X^r,Y^r)$ at the boundary $\{x=0\}\cup\{y=0\}$.
\blem $\P(L^r_{T_0}<\infty)=1$. \elem
\prf Fix $M>0$ and denote $T_M=\inf\{t\ge 0;\; V_t=M\}$. Note that the process $(U_t,V_t)_{t\le T_0}$ is a Brownian motion with oblique reflection in $\{-v\le u\le v\}$, stopped when it hits $(0,0)$.
For $t < T_0$, one has that
\beqarr U_t &=& U_0+W^1_t +\sqrt{2} \int_0^t (1_{\{U_s=-V_s\}}-1_{\{U_s=-V_s\}}) dL^r_s  \\
V_t &=& V_0 + W^2_t,
\eeqarr
where $W^1=\frac{B^2+B^1}{\sqrt{2}}$ and $W^2=\frac{B^2-B^1}{\sqrt{2}}$ are two independent Brownian motions.
Let $h(u,v)=\frac{u^2+v^2}{2v}$. Then, when $u=\pm v$, $\partial_v h(u,v)=0$ and $\partial_u h(u,v)=\pm 1$.
Note also that in $\{-v< u < v\}$,
$$\Delta h (u,v) = \frac{u^2+v^2}{v^3} \le \frac{2}{v}.$$
Applying It\^o's formula, we get for $t<T_0\wedge T_M$ that
$$h(U_t,V_t)\le h(U_0,V_0) + \int_0^t\frac{ds}{V_s} - \sqrt{2} L^r_t + M_t$$
where $M_t$ is a local martingale with quadratic variation given by
$$\langle M\rangle_t = \frac{1}{4}\int_0^t\left(1+\frac{U_s^2}{V_s^2}\right)^2ds .$$
Since this quadratic variation is dominated by $t$, $(M_t)_{t\le T_0}$ is a martingale. This implies that
$$\E\left[h((U,V)_{t\wedge T_0\wedge T_M})\right] \le \E[h(U_0,V_0)]
 + \E\left[\int_0^{t\wedge T_0\wedge T_M}\frac{ds}{V_s}\right] -  \sqrt{2}\E[L_{t\wedge T_0\wedge T_M}].$$
It is well known that $\E\left[\int_0^{T_0\wedge T_M}\frac{ds}{V_s}\right] < \infty$ (for example, 
using It\^o's formula we get for $t<T_0\wedge T_M$,
$$\left. d(V_t\log(V_t))=(1+\log(V_t))dV_t + \frac{dt}{2V_t}.\right)$$

Taking the limit as $t\to\infty$, using dominated and monotone convergence theorems (using that $h(u,v)$ is dominated by $M$ on $\{-v\le u\le v \le M\}$), we get
$$\E\left[h((U,V)_{T_0\wedge T_M})\right] \le \E[h(U_0,V_0)]  + \E\left[\int_0^{T_0\wedge T_M}\frac{ds}{V_s}\right]-\sqrt{2}\E[L^r_{T_0\wedge T_M}].$$
We then must have that $\E[L^r_{T_0\wedge T_M}]<\infty$. So a.s., $L^r_{T_0\wedge T_M}<\infty$ for all $M>0$. This proves the lemma. \qed

\medskip
\blem
Consider again $(X_t,Y_t)$ a solution of the martingale problem \eqref{MPC} for $n=2$. Then $T=\inf\{t\ge 0;\;  X_t=Y_t\}$ is finite a.s.
\elem
\prf Suppose $Y_0>X_0$, so that $(X^r_0,Y^r_0)\in D^+$. To simplify a little, we will also assume that $(X_0,Y_0)\in D$.
Note first that when $T<\infty$, we must have $X_T=Y_T=0$.
Since we are only interested to $(X,Y)$ up to $T$, when $T<\infty$ we will replace  $(X_{T+t},Y_{T+t})_{t>0}$ by $(B_t,B_t)_{t>0}$, with $B$ a Brownian motion independent of $(X_t,Y_t)_{\{t\le T\}}$.
Fix $\eps>0$ and define the sequences of stopping times $\sigma^\eps_k$ and $\tau_k^\eps$, by 
$\tau_0^\eps=0$ and for $k\ge 1$ 
\begin{eqnarray*}
\sigma_{k}^\eps &=& \inf\{t\ge \tau_{k-1}^\eps,\; X_t=\eps \hbox{ or } Y_t=-\eps\}\\
\tau_k^\eps &=& \inf\{t\ge \sigma _{k}^\eps,\; X_t=0 \hbox{ or } Y_t=0\}.
\end{eqnarray*}
Then all these stopping times are finite a.s.
Let us also remark that  $T\not\in\cup_k [\sigma_k^\eps,\tau_{k}^\eps]$  and that $\sum_k (\tau_{k}^\eps\wedge T - \sigma_k^\eps\wedge T)$ converges a.s. towards $T-T_0$ (i.e. the time spent by $(X,Y)$ in $D^c$ before $T$). Take $\alpha>0$. 
Denoting $Z=(X,Y)$, $Z_{\tau_k^\eps}$ is a Markov chain. Note that $N_\eps=\inf\{k,\; Z_{\tau_k^\eps}=(0,0)\}$ is a stopping time for this Markov chain. 
Moreover it is clear that $(\tau_k^\eps-\sigma_k^\eps)_{k\ge 1}$ is a sequence of independent variables, and independent of the Markov chain $Z_{\tau_k^\eps}$. 
And we have for all $k\ge 1$, $\E[e^{-\alpha(\tau_k^\eps-\sigma_k^\eps)}]=e^{-\eps\sqrt{2\alpha}}$. Hence 
\beqarr
\E[e^{-\sum_{k} \alpha(\tau_k^\eps\wedge T - \sigma_k^\eps\wedge T)}]
&=& \E[e^{-\sum_{k= 1}^{N^\eps-1} \alpha(\tau_k^\eps - \sigma_k^\eps)}]\\
&=& \sum_n \P[N_\eps=n] e^{-(n-1)\eps\sqrt{2\alpha}}
\eeqarr
which implies that
$$\E[e^{-\alpha\sum_{k} (\tau_k^\eps\wedge T - \sigma_k^\eps\wedge T)}]=\E[e^{-(N_\eps-1)\eps\sqrt{2\alpha}}].$$

By taking the limit as $\eps\to 0$ in this equality, since $\eps N_\eps$ converges in probability towards $L_T(X)+L_T(Y)$. Denote $L_{T}:=\frac{L_T(X)+L_T(Y)}{2}$. Note that $L_T=L^r_{T_0}$, so that we have
$$\E[e^{-\alpha (T-T_0)}]=\E[e^{-2\sqrt{2\alpha}L^r_{T_0}}].$$
Since $L^r_{T_0}$ is finite a.s., $T-T_0$ is also finite a.s.  \qed

\bprop There exists a unique consistent family $(\P^{(n)}_t,\;n\ge 1)$ of Feller semigroups on $\RR$ such that if $X^{(n)}$ is its associated $n$-point motion started from $x\in\RR^n$
and $T^{(n)}= \inf\{t \ge 0;\;X^{(n)}_t \in \Delta_n\}$, then:
\bdes \iti $(X^{(n)}_t)_{t\le T^{(n)}}$ is distributed like $\PP^{(n),0}_x$.
\itii For $t \ge T^{(n)}$, $X^{(n)}_t\in\Delta_n$.
\edes
Moreover, this family is associated to a coalescing SFM. \eprop
\prf This proposition  follows from Theorem 4.1 in \cite{LJR2004a}. To prove the Feller property, it suffices to check condition (C) of Theorem 4.1 in \cite{LJR2004a}.
Denoting by $(X,Y)$ the two-point motion, condition (C) is verified as soon as for all positive $t$ and $\eps$ (denoting $d(x,y)=|y-x|$)
\begin{equation} \label{cond(C)} \lim_{d(x,y)\to 0}\PP^{(2)}_{(x,y)}[d(X_t,Y_t)>\eps] = 0. \end{equation}
Assume $0<y-x=\alpha<\eps$. Then
\beqarr
\PP^{(2)}_{(x,y)}[d(X_t,Y_t)>\eps]
&=& \PP^{(2)}_{(x,y)}[Y_t-X_t>\eps\hbox{ and } t<T]\\
&\le& \PP^{(2)}_{(x,y)}\left[\sup_{t\le T_0}(Y_t^r-X_t^r)>\eps\right]\\
&\le&   \PP^{(2)}_{(x,y)}\left[\sup_{t\le T_0}V_t>\frac{\eps}{\sqrt{2}}\right].
\eeqarr
This last probability is equal to the probability that a Brownian motion started at $\alpha/\sqrt{2}$ hits $\eps/\sqrt{2}$ before hitting $0$, which is equal to $\alpha/\eps$. This implies \eqref{cond(C)}. \qed

\bprop \label{propsde}Let $\p$ be the SFM associated to $\P^{(n)}$. Then there exist $W^+$ and $W^-$ two independent white noises, $\sigma(\p)$-measurable, such that \eqref{edsCpm} is satisfied.
\eprop
\prf
Define for $s<t$, $W^{\pm}_{s,t}=\lim_{x\to \pm \infty} (\p_{s,t}(x)-x)$. This limit exists a.s. since one can check (using the martingale problem) that for $y>x>0$ (resp. $y<x<0$)  and $t\le \tau^x_s=\inf\{u>s;\; \p_{s,u}(x)=0\}$, $\p_{s,t}(x)-x=\p_{s,t}(y)-y$.
Using that $\langle W^\pm_{s,\cdot},\p_{s,\cdot}(y)\rangle_t = \lim_{x\to\pm\infty} \langle \p_{s,\cdot}(x),\p_{s,\cdot}(y)\rangle_t = 1_{\{\pm \p_{s,t}(y)>0\}}$, 
it is then easy to check \eqref{edsCpm}. \qed

\subsection{Uniqueness}
Let $\wt{\P}^{(n)}$ be another consistent family solving the martingale problem \eqref{MPC}. By uniqueness of the solution of this martingale problem up to the first coalescing time, if  this consistent family is different to the family $\P^{(n)}_t$, then the $n$-point motion should leave $\Delta_n$. After a time change, one may assume it spends no time in $\Delta_n$.
Denote by $(\wt{X},\wt{Y})$ the associated two-point motion starting from $(0,0)$.
Denote $\wt{A}_t=\int_0^t1_{\{(\wt{X}_s,\wt{Y}_s)\in D\}}ds$, $\wt{\kappa}_t= \inf\{s\ge 0;\;\wt{A}_s\ge t\}$ and
$(\wt{X}^r_t,\wt{Y}^r_t)=(\wt{X},\wt{Y})_{\kappa_t}$.
Denote $L_t(\wt{X})$ and $L_t(\wt{Y})$ the local times at $0$ of $\wt{X}$ and of $\wt{Y}$. Define also  $\wt{L}_t$ by $\frac{1}{2} \big( L_t(\wt{X})+L_t(\wt{Y}) \big)$ and $\wt{L}^r_t$ by $\wt{L}^r_{\wt{\kappa}_t}$.
Denote $\wt{U}^r=\frac{\wt{Y}^r+\wt{X}^r}{\sqrt{2}}$ and $\wt{V}^r=\left|\frac{\wt{Y}^r-\wt{X}^r}{\sqrt{2}}\right|$.

\blem
For all $a>0$,  $\wt T_a=\inf\{t>0;\; d(\wt{X}_t,\wt{Y}_t)=a\}$ is infinite a.s.
\elem
\prf Let $\wt T_a^r=\inf\{t>0;\; \wt{V}^r_t)=\frac{a}{\sqrt{2}}\}$. Then one has $\wt A_{\wt T_a} = \wt T_a^r$.
Let $\eps<\frac{a}{\sqrt{2}}$.
Define $\sigma_n^\eps$ and $\tau_n^\eps$ by $\tau_0^\eps=0$ and for $n\ge 1$,
\begin{eqnarray*}
\sigma_n^\eps &=& \inf\{t\ge \tau_{n-1}^\eps;\; \wt{V}^r_t=\eps\}\\
\tau_n^\eps &=& \inf\{t\ge \sigma_n^\eps;\; \wt{V}^r_t=0\}.
\end{eqnarray*}
Recall $h(u,v)=\frac{u^2+v^2}{2v}$
Applying It\^o's formula on the time intervals $[\sigma_n^\eps\wedge \wt T_a^r,\tau_n^\eps\wedge \wt T_a^r]$, denoting $\wt{Z}^r=(\wt{U}^r,\wt{V}^r)$, and using $\Delta h(u,v)\ge \frac{1}{v}$ for the second term, 
\begin{eqnarray}
h\big(\wt{Z}^r_{\wt T^r_a}\big)-h(0)
&\ge& \sum_{n\ge 1} \left(h\big(\wt{Z}^r_{\tau_n^\eps\wedge \wt T_a^r}\big)-h\big(\wt{Z}^r_{\sigma_{n-1}^\eps\wedge \wt T_a^r}\big)\right)\label{term1}\\
&+& \sum_{n\ge 1} \int_{\sigma_n^\eps\wedge \wt T_a^r}^{\tau_n^\eps\wedge \wt T_a^r} \frac{1}{2\wt{V}^r_s} ds \label{term2}\\
&-& \sum_{n\ge 1} \left(\wt{L}^r_{\tau_n^\eps\wedge \wt T_a}-\wt{L}^r_{\sigma_n^\eps\wedge \wt T_a}\right)\label{term3}\\
&+& \sum_{n\ge 1} \left(\wt{M}^{n,\eps}_{\tau_n^\eps\wedge \wt T_a^r}-\wt{M}^{n,\eps}_{\sigma_n^\eps\wedge \wt T_a^r}\right)\label{term4}
\end{eqnarray}
where for all $n\ge 1$ and all $\eps>0$, $\wt{M}^{n,\eps}$ is a martingale whose quadratic variation  is such that $\frac{d}{dt}\langle \wt{M}^{n,\eps}\rangle_t\le 1$. 
Since  $\int_0^{T^B_r}\frac{ds}{|B_s|}=\infty$ a.s. (see Corollary 6.28 Chapter 3 in \cite{KS}), with $T^B_r$ the first time a Brownian motion $B$ hits $-a/\sqrt{2}$ or $a/\sqrt{2}$, we get that
term \eqref{term2} goes to $\infty$ as $\eps\to 0$. Terms \eqref{term1} is positive and is dominated by $\eps \times \#\{n;\;\tau_n^\eps \le \wt T_a^r\}$ which converges in probability towards $L^0_{\wt T_a^r}$, with $L^0$ the local time at $0$ of $\wt{V}$, which is a reflected Brownian motion. Denote by $M^\eps$ the martingale defined by
$$M^\eps_t  = \sum_{n\ge 1} (\wt{M}^{n,\eps}_{\tau_n^\eps\wedge t}-\wt{M}^{n,\eps}_{\sigma_n^\eps\wedge t}).$$
Then if $\eps'<\eps$, $\langle M^\eps - M^{\eps'}\rangle_t=\langle M^{\eps'}\rangle_t-\langle M^\eps\rangle_t$. 
Using that $\frac{d}{dt}\langle \wt{M}^{n,\eps}\rangle_t\le 1$,  we get $\langle M^\eps - M^{\eps'}\rangle_{\wt T_a^r}\le \sum_{n\ge 1} (\sigma_n^\eps\wedge \wt T_a^r - \tau_{n-1}^\eps\wedge \wt T_a^r)$.
Since $\E[\sum_{n\ge 1} (\sigma_n^\eps\wedge \wt T_a^r - \tau_{n-1}^\eps\wedge \wt T_a^r)]$ converges towards $0$ as $\eps\to 0$, we get that $M^\eps_{\wt T_a^r}$ is a Cauchy sequence in $L^2\big(\wt{\PP}^{(2)}_{(0,0)}\big)$ and thus \eqref{term4}, which is equal to $M^\eps_{\wt T_a^r}$, converges in $L^2\big(\wt{\PP}^{(2)}_{(0,0)}\big)$ .
This implies that term \eqref{term3} converges in probability towards $\infty$. 
Since it also converges towards $\wt{L}^r_{\wt T_a^r}$. One concludes that $\wt{L}^r_{\wt T_a^r}=\infty$ a.s. 
Since $\wt T_a=\wt{\kappa}_{\wt T_a^r}$, $\wt{L}^r_{\wt T_a^r}=\wt{L}_{\wt T_a}$ and therefore $\wt{L}_{\wt T_a}=\infty$ a.s.
Since $\wt{X}$ and $\wt{Y}$ are two Brownian motions, one must have $\wt T_a=\infty$ a.s. \qed

\medskip
This lemma implies that after hitting $\Delta_2=\{x=y\}$, the two-point motion stays in $\Delta_2$. Thus $\wt{\P}^{(2)}=\P^{(2)}$. The same argument applies to prove that $\wt{\P}^{(n)}=\P^{(n)}$. This proves that there exists only one flow which is a SFK solution of the $(\frac{1}{2}\Delta,C_\pm)$-SDE. Moreover this solution is a SFM  $\p^\pm$.

To prove $\p^\pm$ is a Wiener solution, 
define $\bar{K}$ the SFK obtained by filtering $\p^\pm$ with respect to the noise $N^{\pm}$   generated by $W^+$ and $W^-$. 
Then, like for Proposition \ref{propsde} we have $\sigma(W^+,W^-)\subset \sigma(\bar{K})$, and applying Proposition \ref{prop:edsfilt}, $\bar{K}$ also solves the $(\frac{1}{2}\Delta,C_\pm)$-SDE driven by $W^+1_{\RR^+}+W^-1_{\RR^-}$. Since there exists a unique solution, $\bar{K}$ is distributed like $\delta_{\p^\pm}$. The noise of $\bar{K}$ being $N^{\pm}$, this implies that $\p^\pm$ is measurable with respect to $\sigma(W^+,W^-)$, i.e. $\p^\pm$ is a Wiener solution. This concludes the proof of Theorem \ref{thm2}.


\section{$C_+(x,y)=1_{x>0}1_{y>0}$.}
\subsection{The Wiener solution obtained by filtering.}\label{wsc+}
Let $\p^\pm$ be the SFM solution of the $(\frac{1}{2}\Delta,C_\pm)$-SDE, and $N^+$ the noise generated by $W^+$. 
Define $K^+$ the SFK obtained by filtering $\p^\pm$ with respect to $N^+$.
Note that for $x\in\RR$, $s\in\RR$ and $t\in [s,\tau_s(x)]$ with $\tau_s(x)=\inf\{u;\; W^+_{s,u}=-x\}$, we have $\p^{\pm}_{s,t}(x)=x+W^+_{s,t}$ and thus that $K^+_{s,t}(x)=\delta_{x+W^+_{s,t}}$. This clearly implies that $W=W^+ 1_{\RR^+}$ is a $\cF^{K^+}$-white noise. 
Then (see section \ref{filtering}) $K^+$ is a Wiener solution of the $(\frac{1}{2}\Delta,C_+)$-SDE driven $W=W^+ 1_{\RR^+}$.
Proposition  \ref{21} implies it is the unique Wiener solution.

\subsection{Construction of the $n$-point motions.} \label{sec:npt0}
Let $D_n=\{x\in\RR^n; \; \exists i\neq j,\; x_i=x_j\ge 0\}$. Let $x\in\RR^n\backslash  D_n$.

The $n$-point motion
has the property that at any time, the points located on the
positive half line will move parallel to $W^+$ and the points
located on the negative half line will follow independent Brownian
motions. Its law is determined up to the first time it hits $D_n$. Indeed,
denote $I_-=\{j;\; x_j\le 0\}$ and $I_+=\{j;\; x_j>0\}$.  Let $i$  be
such that $x_i=\max\{x_j;\; j\in I_-\}$ when $I_-\neq \emptyset$ and
$x_i=\min\{x_j;\; 1\le j\le n\}$ when $I_-=\emptyset$. Let
$(B^1,\dots,B^n,B)$  be $n+1$ independent Brownian motions. Define
for $t>0$, the processes \beqarr
{X}^0_j(t)&=& x_j+B^j(t) \qquad\hbox{ for } \quad j\in I_-\cup \{i\}\\
{B}^{0,+}_t &=& \int_0^t 1_{\{{X}^0_i(s)<0\}}dB_s  + \int_0^t 1_{\{{X}^0_i(s)>0\}}dB^i_s\\
{X}^0_j(t) &=&  x_j+{B}^{0,+}_t  \qquad \hbox{ for } \quad j\in I_+\backslash\{i\}
\eeqarr
Set
 $$\tau_1=\inf\{t>0;\; \exists j\neq i,\; {X}^0_{j}(t)=0 \}$$ and set for $t\leq \tau_1$, ${X}_t={X}^0_t$ and ${B}^+_t={B}^{0,+}_t$.

Assume now that $(\tau_k)_{k\le \ell}$ and $({X}(t),{W}^+(t))_{t\le\tau_\ell}$ have been defined such that a.s.
\begin{itemize}
\item $(\tau_k)_{1\le k\le\ell}$ is an increasing sequence of stopping times with respect to the filtration associated to ${X}$;
\item ${X}_{\tau_k}\not\in D_n$ for $1\le k\le \ell$;
\item for all $k\le \ell$, there exists an integer $i_k$ such that ${X}_{i_k}(\tau_k)=0$.
\end{itemize}
We then define $({X}_t,{B}^+_t)_{\tau_\ell<t\le \tau_{\ell+1}}$ as is defined $({X}_t,B^+_t)_{0<t\le \tau_{1}}$ by replacing $i$ by $i_\ell$, $x$ by ${X}_{\tau_\ell}$ and $(B^1_\cdot,\dots,B^n_\cdot,B_\cdot)$ by $(B^1_{\tau_\ell+\cdot},\dots, B^n_{\tau_\ell+\cdot},B_{\tau_\ell+\cdot})-(B^1_{\tau_\ell},\dots, B^n_{\tau_\ell},B_{\tau_\ell})$.
Let ${T}=\lim_{\ell\to\infty}\tau_\ell$. Denote by $\PP^{(n),0}_x$ the law of $(X_t)_{t\leq T}$.

\blem \label{lem:MPunicity2} Let $X^{(n)}$ be a solution to the martingale problem \eqref{MPC}, with $X^{(n)}_0=x$. Let $T^{(n)}=\inf\{t;\;\exists i\neq j,\; X^{(n)}_t\in D_n\}$. Then  $(X^{(n)}_t)_{t\le T^{(n)}}$ is distributed like $\P^{(n),0}_x$. \elem
\prf Let $x=(x_1,\dots,x_n)\in\RR^n\backslash D_n$.
As before $I_-=\{j;\; x_j\le 0\}$ and $I_+=\{j;\; x_j>0\}$. Let $i$  be such that $x_i=\max\{x_j;\; j\in I_-\}$ when $I_-\neq \emptyset$ and $x_i=\min\{x_j;\; j\in I\}$ when $I_-=\emptyset$. Let $X^{(n)}$ be a solution of the martingale problem.
This implies that for all $j$, $X^{(n)}_j$ is a Brownian motion and for all $j\neq k$,
\begin{equation}
\label{crochet2}\langle X^{(n)}_j,X^{(n)}_k\rangle_t = \int_0^t 1_{\{X^{(n)}_j(s)>0\}}1_{\{X^{(n)}_k(s)>0\}} ds.\end{equation}

Let $(B^{0,1},\dots,B^{0,n},B^0)$  be $n+1$ independent Brownian motions, and independent of $X^{(n)}$.
Define $(B^1,\dots,B^n)$ by
$$B^j_t= \int_0^t 1_{\{X^{(n)}_j(s)<0\}} dX^{(n)}_{j}(s) + \int_0^t 1_{\{X^{(n)}_j(s)>0\}}dB^{0,j}_s.$$
Note that $(B^1,\dots,B^n)$ are $n$ independent Brownian motions.

Set $$\tau_1=\inf\{t>0;\; \exists j\neq i,\; X^{(n)}_{j}(t)=0\}.$$ Define for $t\le \tau_1$,
$B_t=B^0_t$ when $I_+=\emptyset$ and
$$B_t = \int_0^t 1_{\{X^{(n)}_i(s)<0\}} dX^{(n)}_k(s) + \int_0^t 1_{\{X^{(n)}_i(s)>0\}}dB^0_s$$
when $I_+\neq \emptyset$ and where $k\in I_+$ (one can choose for example $k$ such that $x_k=\max\{x_j;\;j\in I_+\}$).
Define also for $t\le \tau_1$
$${B}^{+}_t = \int_0^t 1_{\{{X}^{(n)}_i(s)<0\}}dB_s  + \int_0^t 1_{\{{X}^{(n)}_i(s)>0\}}dB^i_s.$$
Note that for $t\le \tau_1$, $\langle B,B^j\rangle_t=0$ for all $j$ and that
\beqarr
{X}^{(n)}_j(t)&=& x_j+B^j_t \qquad\hbox{ for } \quad j\in I_-\cup \{i\}\\
{X}^{(n)}_j(t) &=& x_j+B^{+}_t \qquad \hbox{ for } \quad j\in I_+\backslash\{i\}.
\eeqarr

Assume now that $(\tau_k)_{k\le \ell}$ and $(B_t,B^+_t)_{t\le\tau_\ell}$ have been defined such that a.s.
\begin{itemize}
\item $(\tau_k)_{1\le k\le\ell}$ is an increasing sequence of stopping times with respect to the filtration associated to ${X^{(n)}}$;
\item $X^{(n)}_{\tau_k}\not\in \Delta_n$ for $1\le k\le \ell$;
\item for all $1\le k\le \ell$, there exists an integer $i_k$ such that $X^{(n)}_{i_k}(\tau_k)=0$.
\end{itemize}
We then define $(B_t,B^+_t)_{t\in ]\tau_\ell, \tau_{\ell+1}]}$ as is defined $(B_t,B^+_t)_{0<t\le \tau_{1}}$ by replacing $i$ by $i_\ell$, $x$ by ${X^{(n)}}_{\tau_\ell}$.
Note that ${T^{(n)}}=\lim_{\ell\to\infty}\tau_\ell$. Define for $t\ge T^{(n)}$, $B_t=B_{T^{(n)}}+B^0_t-B^0_{T^{(n)}}$. Then $(B^1,\dots,B^n,B)$ are $n$ independent Brownian motions.
We finally remark that $(X^{(n)}_t)_{t\le T^{(n)}}$ is defined $\PP^{(n),0}_x$. This proves the lemma.
\qed

\subsection{Proof of Theorem \ref{thm4}}
Denote by $\P^{(n)}$ the family of consistent Feller semigroups associated to $K^+$.
To this family of semigroups, we associate a unique consistent family of coalescing $n$-point motions,  $\P^{(n),c}$ (see Theorem 4.1 in \cite{LJR2004a}), with the property that the law of this $n$-point motion before hitting $\Delta_n$ is described by Lemma \ref{lem:MPunicity2}. These semigroups are Fellerian and since the two-point motion hits $\Delta_2=\{x=y\}$, the family of semigroups are associated to a coalescing SFM $\p^+$.

Proof of (i): follow the proof of Proposition \ref{propsde}.

Proof of (ii): It is a consequence of Theorem 4.2 in \cite{LJR2004a}.
 \qed

\subsection{Proof of Theorem \ref{thm1}}
The existence (and uniqueness) of a Wiener solution $K^+$ was proved in section
\ref{wsc+}. By construction, (iii) holds. 
The SFK $K^+$ is diffusive since the two-point
motion clearly leaves $\Delta_2$ (it behaves as a Brownian motion in
$D=\{x<0\hbox{ or }y<0\}$). So (ii) is proved. To finish the proof of (i), it remains to show that $K^+$ is the unique solution. Let $K$ be another solution. The flow obtained by filtering $K$ with respect to $W^+$ is a Wiener solution. Since there exists at most one Wiener solution, this flow is distributed like $K^+$. For simplicity we also denote this flow $K^+$.
Denote by $\P^{(2)}$ (resp. $\P^{(2,+)}$) the semigroup of the two-point of $K$ (resp. of $K^+$). Note that if  $\P^{(2)}=\P^{(2,+)}$, then $K=K^+$. Indeed,
\begin{eqnarray*}
\E[(K_{s,t}f(x)-K^+_{s,t}f(x))^2]
&=& \P^{(2)}_{t-s}(f\otimes f) (x,x) -2\E[K_{s,t}f(x)K^+_{s,t}f(x)] \\
&&+  \P^{(2,+)}_{t-s}(f\otimes f) (x,x)\\
&=& \P^{(2)}_{t-s}(f\otimes f) (x,x) -2\E[K^+_{s,t}f(x)K^+_{s,t}f(x)]\\
&& +  \P^{(2,+)}_{t-s}(f\otimes f) (x,x)\\
&=& \P^{(2)}_{t-s}(f\otimes f) (x,x) -  \P^{(2,+)}_{t-s}(f\otimes f) (x,x)\\
&=& 0.
\end{eqnarray*}

Denote
\begin{eqnarray*}
D^r &=& \{(x,y)\in\RR^2;\; x\le 0 \hbox{ or } y\le 0\}\\
D^+ &=& \{(x,y)\in\RR^2;\; x\ge 0 \hbox{ and } y\ge 0\}\\
B &=& \{(x,y)\in\RR^2;\; x=0 \hbox{ and } y\ge 0\}\\
&& \cup \{(x,y)\in\RR^2;\; y=0 \hbox{ and } x\ge 0\}
\end{eqnarray*}
so that $B$ is the boundary of $D^r$ and $D^+$, and $\RR^2=D^r\cup B\cup D^+$.

Let $Z=(X,Y)$ be a solution to the martingale problem:
\begin{equation} \label{MPC+}
f(Z_t)-\int_0^t \frac{1}{2}\Delta f (Z_s) ds - \int_0^t 1_{\{Z_s\in D^+\}} \frac{\partial^2}{\partial_x\partial_y}f(Z_s) ds \end{equation}
is a martingale for all $f\in C^2_0(\RR^2)$.
For $\eps>0$, let
\begin{eqnarray*}
B^\eps &=& \{z\in D^+;\; d(z,B)=\eps\} \end{eqnarray*}
with $d$ the Euclidean distance on $\RR^2$ (note that when $z=(x,y)\in D^+$, $d(z,B)=x\wedge y$).

Define $A^r_t=\int_0^t 1_{\{Z_s\in D^r\}} ds$ and $A^+_t=\int_0^t 1_{\{Z_s\in D^+\}} ds$, the amount of time $Z$ spends in $D^r$ and $D^+$ up to time $t$. Denote by $\kappa^r_t$ and by $\kappa^+_t$ the inverses of $A^r_t$ and of $A^+_t$:
\begin{eqnarray*}
\kappa^r_t &=& \inf\{s\ge 0;\; A^r_s>t\}\\
\kappa^+_t &=& \inf\{s\ge 0;\; A^+_s>t\}
\end{eqnarray*}
Denote by $Z^r$ and by $Z^{+}$ the processes in $D^r$ and $D^+$ defined by
$$ Z^r_t=Z_{\kappa^r_t}\quad \hbox{ and }\quad Z^{+}_t=Z_{\kappa^+_t}. $$

Define the sequences of stopping times $\sigma_k^\eps$ and $\tau_k^\eps$ by $\sigma_0^\eps=0$ and for $k\ge 0$
\begin{eqnarray*}
\tau_k^\eps &=& \inf\{t\ge \sigma_k^\eps;\; Z_t\in B\}\\
\sigma_{k+1}^\eps &=& \inf\{t\ge \tau_k^\eps;\; Z_t\in B^\eps\}.
\end{eqnarray*}

It is easy to see that for $k\ge 1$, in the time interval $[\sigma_k^\eps,\tau_k^\eps]$, $Z_t\in D^+$, $Y_t-X_t$ remains constant and that $d(Z,B)$ is a Brownian motion stopped when it hits $0$ (note that $d(Z,B)=X\wedge Y$).
Thus, if $R_t = d(Z^{+}_t,B)$, then $(R_t)_{t\le A^+_\infty}$ is a Brownian motion instantaneously reflected at $0$ (since $R_t-R_0=\lim_{\eps\to 0} \sum_{k\ge 0} \left(d(Z_{\sigma_k^\eps\wedge \kappa^+_t},B) - d(Z_{\tau_k^\eps\wedge \kappa^+_t},B)\right)$).

Denote now by $L^+_t$ the local time at $0$ of $R_t$. Then, for any random times $T$, $L^+_{A^+_T}$ is the limit in probability as $\eps\to 0$ of $\eps\times N_T^\eps$ where
$$N^\eps_T=max\{k;\; \sigma_k^\eps <T\}.$$
Denote also $L^r_t=L^+_{A^+_{\kappa^r_t}}$.

\blem For all $f\in C^2_0(\RR^2)$ (resp. $f\in C^2(\RR^2)$), we have
\begin{equation}
f(Z^r_t)-\frac{1}{2}\int_0^t \Delta f(Z^r_s) ds - \int_0^t (\partial_x+\partial_y)f(Z^r_s)dL^r_s
\end{equation}
is a martingale (resp. a local martingale).
\elem
\prf Take $f\in C^2_0(\RR^2)$ and assume $Z_0\in D^+$ (we leave to the reader the case $Z_0\not\in D^+$). Then 
\begin{eqnarray}
f(Z^r_t)-f(Z^r_0)
&=& \sum_{k\ge 1} \left(f(Z_{\tau_k^\eps\wedge \kappa^r_t}) - f(Z_{\sigma_k^\eps\wedge \kappa^r_t})\right)\label{*1}\\
&+& \sum_{k\ge 0} \left(f(Z_{\sigma_{k+1}^\eps\wedge \kappa^r_t}) - f(Z_{\tau_k^\eps\wedge \kappa^r_t})\right)\label{*2}
\end{eqnarray}
Note that for $k\ge 1$ and $\kappa^r_t>\sigma_k^\eps$ (which implies $\kappa^r_t>\tau_k^\eps$), $Z_{\tau_k^\eps\wedge \kappa^r_t} - Z_{\sigma_k^\eps\wedge \kappa^r_t} = -\eps v$, with $v=(1,1)$.
Using Taylor expansion, we have that the first term \eqref{*1} is equal to
$$-\eps \sum_{k=1}^{N^\eps_{\kappa^r_t}} (\partial_x+\partial_y)f(Z_{\sigma_k^\eps}) + O(\eps^2\times N^\eps_{\kappa^r_t}).$$
More precisely, $O(\eps^2\times N^\eps_{\kappa^r_t})\le \frac{\eps^2 N^\eps_{\kappa^r_t}}{2}\|f\|_{2,\infty}$. This implies that $O(\eps^2\times N^\eps_{\kappa^r_t})$ converges in probability towards $0$ as $\eps \to 0$ (recall $\eps\times N^\eps_{\kappa^r_t}$ converges towards $L^r_t=L^+_{A^+_{\kappa^r_t}}$).
Thus the first term converges towards
$$-\int_0^t(\partial_x+\partial_y)f(Z^r_s)dL^r_s.$$
(it is obvious if $(\partial_x+\partial_y)f(z)=1$ when $z\in B$).
Indeed:
we claim that for $H$ a bounded continuous process,
$\eps \sum_{k=1}^{N^\eps_{t}} H_{\sigma_k^\eps}$
converges in probability towards $\int_0^{A^+_t} H_s dL_s$. This holds for $H=\sum_i H_{t_i}1_{]t_{i-1},t_i]}$. Every bounded continuous process can be approached by a sequence of processes $H^n=\sum_i H_{i2^{-n}}1_{](i-1)2^{-n},i2^{-n}]}$. Thus we prove the claim by density. Then we apply this result by taking $H_s=(\partial_x+\partial_y)f(Z_s)$ to prove the convergence of the first term \eqref{*1}.

Denote by $M^f$ the martingale \eqref{MPC+}. Then the second term \eqref{*2} is equal to
\begin{eqnarray}
&& \sum_{k\ge 1} (M^f_{\sigma_k^\eps\wedge \kappa^r_t} - M^f_{\tau_{k-1}^\eps\wedge \kappa^r_t}) \label{**1}\\
&+& \frac{1}{2} \sum_{k\ge 1}\int^{\sigma_k^\eps\wedge \kappa^r_t}_{\tau_{k-1}^\eps\wedge \kappa^r_t}
\Delta f(Z_s) 1_{\{Z_s\in D\}} ds \label{**2}\\
&+&  \sum_{k\ge 1}\int^{\sigma_k^\eps\wedge \kappa^r_t}_{\tau_{k-1}^\eps\wedge \kappa^r_t}
\frac{\partial^2}{\partial_x\partial_y}f(Z_s) 1_{\{Z_s\in D^+\}} ds\label{**3}
\end{eqnarray}
Note that
$$\eqref{**2} = \frac{1}{2} \int_0^{t} \Delta f(Z^r_s) ds.$$

Note that for a constant $C$ depending only on $f$, and denoting $D_\eps^+=\{z\in D^+;\; d(z,B)\le \eps\}$,
$$|\eqref{**3}|\le C \int_0^{\kappa^r_t} 1_{\{Z_s\in D^+_\eps\}}ds$$
which converges a.s. towards $0$ as $\eps\to 0$.

Denote by $M^\eps_t$ the term \eqref{**1}. It is a martingale for all $\eps>0$ and one can check that there exists a constant $C<\infty$ depending only on $f$ such that
\begin{eqnarray*}
\frac{d}{dt} \langle M^f\rangle_t  &=&   1_{\{Z_t\in D^r\}} |\nabla f(Z_t)|^2 + 1_{\{Z_t\in D^+\}} |(\partial_x+\partial_y) f(Z_t)|^2 \\
&\le& C
\end{eqnarray*}
Moreover, for $\eps'<\eps$,
\begin{eqnarray*}
\langle M^\eps-M^{\eps'}\rangle_t
&=& \langle M^\eps\rangle_t - \langle M^{\eps'}\rangle_t \\
&=& \sum_{k\ge 0} \int_{\tau_k^{\eps}\wedge\kappa^r_t}^{\sigma_{k+1}^{\eps}\wedge\kappa^r_t} \frac{d}{dt} \langle M^f\rangle_t - \sum_{k'\ge 0} \int_{\tau_{k'}^{\eps'}\wedge\kappa^r_t}^{\sigma_{k'+1}^{\eps'}\wedge\kappa^r_t} \frac{d}{dt} \langle M^f\rangle_t.
\end{eqnarray*}
Since $\{s\ge 0;\; Z_s\in D^r\cup D^+_\eps\} \supset \cup_{k\ge 0} [\tau_k^{\eps},\sigma_{k+1}^{\eps}] \supset \cup_{k'\ge 0} [\tau_{k'}^{\eps'},\sigma_{k'+1}^{\eps'}] \supset \{s\ge 0;\; Z_s\in D^r\}$, we get that
\begin{eqnarray*}
\langle M^\eps-M^{\eps'}\rangle_t
&\le& C \int_0^{\kappa^r_t} 1_{\{Z_s\in D^+_\eps\}}ds.
\end{eqnarray*}
Since this converges a.s. towards $0$ as $\eps\to 0$, uniformly in $\eps'<\eps$, we have that $M^\eps$ is a Cauchy sequence in the space of square integrable martingales. Thus $M^\eps$ converges towards a martingale $M$. This proves the Lemma. \qed

\medskip
This Lemma implies that $Z^r$ is a Brownian motion with oblique reflection in the wedge $D^r$ (see \cite{VW85}, it corresponds to the case $\xi=3\pi/2$, $\theta_1=\theta_2=\pi/4$ and $\alpha=1/3$): $\int_0^t 1_{\{Z^r_s=0\}}ds = 0$ and $Z^r$ is a solution to the sub-martingale problem:
\begin{equation}
f(Z^r_t)-\frac{1}{2}\int_0^t \Delta f(Z^r_s) ds
\end{equation}
is a sub-martingale for all $f$ constant in the neighborhood of $0$ and $f\in C^2_0(D)$ such that $(\partial_x+\partial_y)f(z)\ge 0$  for $z\in D$.

In \cite{VW85}, it is proved that for all initial value $Z^r_0$, there is a unique solution to this sub-martingale problem. Thus the law of $Z^r$ is uniquely determined by $Z^r_0$.

Applying also the Lemma to the function $f(x,y)=x$, we see that
$X^r_t+L^r_t$ is a local martingale (it is actually a true
martingale). This implies that $X^r$ is a semimartingale and that in
the Doob-Meyer decomposition of $X^r$, $-L^r$ is its compensator.
Thus $L^r$ can be recovered from $X^r$. Note that this gives a proof
that $Z^r$ is a semimartingale (see also \cite{W85}).

Note that $L^r_{A^r_t}=L^+_{A^+_t}$. Thus, if $A^+_\infty=\infty$ then $L^+_{A^+_\infty}=\infty$ (since $\{0\}$ is recurrent for $R$) and $L^r_{A^r_\infty}=\infty$. Therefore $A^+_\infty=\infty$ implies that $A^r_\infty$.
On the converse, if $A^r_\infty=\infty$, then $L^r_\infty=\infty$ (since $B$ is recurrent for $Z^r$) and $L^+_{A^+_\infty}=\infty$. Therefore $A^r_\infty=\infty$ implies that $A^+_\infty=\infty$. Since $A^r_\infty+A^+_\infty=\infty$, we must have $A^r_\infty=A^+_\infty=\infty$.
It is easy to see that the processes $Z^r$ and $R$ are independent. We have that
$$Z_t=Z^r_{A^r_t}+R_{A^+_t}v$$
with $v=(1,1)$.

Set $L_t=L^r_{A^r_t}(=L^+_{A^+_t})$.  Define $T_\ell$, $T^r_\ell$ and $T^+_\ell$ by
\begin{eqnarray*}
T_\ell &=& \inf\{t\ge 0;\;L_t>\ell\} \\
T^r_\ell &=& \inf\{t\ge 0;\;L^r_t>\ell\} \\
T^+_\ell &=& \inf\{t\ge 0;\;L^+_t>\ell\}
\end{eqnarray*}
Then $T_\ell=A^r_{T_\ell}+A^+_{T_\ell}=T^r_\ell+T^+_\ell$. Thus the
process $L_t$ is $\sigma(Z^r,R)$-measurable since $L_t =
\inf\{\ell;\;T_\ell>t\} =  \inf\{\ell;\;T^r_\ell+T^+_\ell>t\}$. Now,
$A^r_t=T^r_{L_t}$ and $A^+_t=T^+_{L_t}$. Since
$$Z_t=Z^r_{A^r_t}+R_{A^+_t}v,$$ we see that $Z$ is
$\sigma(Z^r,R)$-measurable. Thus, the law of $Z$ is uniquely
determined. This proves that $\P^{(2)}=\P^{(2,+)}$. \qed

\bibliographystyle{plain}
\bibliography{edspm}

\end{document}